\documentclass[10pt]{article}
\usepackage[all]{xy}
\usepackage{amsmath}
\usepackage{amsfonts}
\usepackage{amssymb}
\usepackage{amscd}
\usepackage{amsthm}
\usepackage{latexsym}
\usepackage{amsbsy}
\usepackage{graphicx}
\usepackage{epstopdf}
\usepackage{seqsplit}
\usepackage{color}
\usepackage{float}
\usepackage{subfig}
\usepackage{authblk}
\usepackage{natbib}

\AppendGraphicsExtensions{.gif}

\newtheorem{lemma}{Lemma}[section]
\newtheorem{proposition}{Proposition}[section]

\newtheorem{example}{Example}[section]

\newtheorem{remark}{Remark}[section]

\newcommand{\R}{\mathbb{R}}

\newcommand{\tX}{\widetilde X}
\newcommand{\tY}{\widetilde Y}
\newcommand{\tp}{\widetilde p}
\newcommand{\tq}{\widetilde q}
\newcommand{\tr}{\widetilde r}
\newcommand{\mr}{\bar{m}}
\newcommand{\al}{\alpha}
\newcommand{\alzeroone}{0\leq \alpha \leq 1}

\newcommand{\bE}{\mathbb{E}}

\providecommand{\keywords}[1]
{
  \small	
  \textbf{\textit{Keywords---}} #1
}

\title{
Fuzzy Gaussian mixture optimization of the newsvendor problem: mixing online reviews and judgemental demand data}

\author[1]{Farzad Fathizadeh}
\author[2]{Jean Savinien}
\author[2]{Yacine Rekik}
\affil[1]{Computational Foundry, Swansea University, UK}
\affil[2]{emlyon Business School, France}

\date{April 2021}

\begin{document}

\maketitle

\begin{abstract}
Motivated by the increasing exposition of decision makers to both statistical and judgemental based sources of demand information, we develop in this paper a fuzzy Gaussian Mixture Model (GMM) for the newsvendor permitting to mix probabilistic inputs with a subjective weight modelled as a fuzzy number. The developed framework can model for instance situations where sales are impacted by customers sensitive to online review feedbacks or expert opinions. It can also model situations where a marketing campaign leads to different stochastic alternatives for the demand with a fuzzy weight. Thanks to a tractable mathematical application of the fuzzy machinery on the newsvendor problem, we derived the optimal ordering strategy taking into account both probabilistic and fuzzy components of the demand. We show that the fuzzy GMM can be rewritten as a classical newsvendor problem with an associated density function involving these stochastic and fuzzy components of the demand. The developed model enables to relax the single modality of the demand distribution usually used in the newsvendor literature and to encode the risk 
attitude of the decision maker.

\end{abstract}


\section{Research motivation and contribution}

The classic newsvendor problem (\cite{Arrowetal1951}) models a situation where an inventory manager must decide upon an order quantity to trade-off overage and underage risks against an unknown demand from a known distribution. This classic problem is simple to solve where the optimal order quantity is derived by the critical fractile of the demand distribution, but it has several variations. There is now a comprehensive list of extensions on such problems. For further reading on the classical newsvendor problem and its main extensions, the reader is referred to \cite{Khouja1999} and \cite{Quinetal2011}, and the numerous references therein.

By having a look at this comprehensive literature, one could observe that the vast majority of newsvendor investigations consider a given probability distribution for the demand. Indeed, all the newsvendor extensions have been mainly performed under the probabilistic framework under a demand characterized by a probability distribution, and an objective function modelled by an expected utility.
In practice however, the uncertainty of the demand for single periods, especially for fashion and high tech products, makes it difficult to obtain accurate forecasts. This is nowadays more true particularly with consumers trending to provide real feedbacks on their purchases and the increasing possibility of influencing the behavior of forthcoming consumers of the sales period.
\cite{Scarf1958} writes: “we may have reason to suspect that the future demand will come from a distribution that differs from that governing past history in an unpredictable way”. More recently, \cite{BertsimasandThiele2014} describe the need for a non-probabilistic theory as “pressing”. These new trends of consumer behaviors strengthen the unpredictability of the demand and provides consequently a strong incentive for the decision maker to deploy alternative solutions, which should take into account both statistical data on demand as well as subjective judgment on consumer behaviors. There is a frequent trend where forecasters and stock managers judgmentally adjust a statistical forecast or a replenishment decision (\cite{kholidasari_review_2018}). \cite{syntetos_effects_2016} reported that judgmentally adjusted replenishment orders may improve inventory performance in terms of reduced inventory investments (costs).

To tackle the demand unpredictability issue, some investigations consider a free distribution approach for the newsvendor problem. \cite{Scarf1958} is the first to give a closed-form solution to the newsvendor problem when only the mean and the variance of the demand  are assumed to be known. The investigation of \cite{GallegoandMoon1993} considered a distribution-free modeling and provided an extension to Scarf’s solution. \cite{yue_expected_2006} assumed that the demand density function belongs to a specific family of density functions. 
However, with the emergence of social media and IT capabilities, decision makers are facing two types of uncertainty: 
\begin{itemize}
\item Statistical data which can result in a probability distribution function.  
\item  Subjective and human judgement based uncertainty related to past customer reviews and expert judgements which may impact the future sales. 
\end{itemize}
The nature of the latter uncertainty could be handled by average scores attributed by consumers and/or linguistic terms to express their satisfaction of the product. Considering the fact that a distribution-free model does not completely solve the modeling of the subjective uncertainty, another set of investigations  introduced the fuzzy theory, as a promising instrument, to describe and treat the uncertainty in these cases.

Several investigations have analyzed newsvendor problems under a fuzzy demand assumption. \cite{petrovic_fuzzy_1996} considered a fuzzy newsvendor model where the overage cost, the shortage cost, and the demand are fuzzy numbers, and the optimal order quantity is obtained by the defuzzification of total costs.

\cite{li_fuzzy_2002} tried to integrate both probabilistic and fuzzy uncertainties by studying two models: in the first one, the demand was assumed probabilistic while the costs were fuzzy, and in the second one, they assumed that the costs were deterministic but the demand was fuzzy. Other investigations considered a fuzzy assumption for the demand in order to extend the newsvendor problem for the case of product substitution  (\cite{dutta_incorporating_2010}), quantity discount newsvendor problem (\cite{chen_analysis_2011}) or to consider a decentralized supply chain (\cite{xu_optimal_2008}, \cite{ryu_fuzzy_2010}). 
It is worthwhile to notice that these newsvendor investigations handled the demand uncertainty by either a probabilistic approach or by a fuzzy approach. Despite the fact that decision makers are rather exposed to both statistical and subjective type demand data, to the best of our knowledge, there is no published research where these two types of demand uncertainty are mixed. Our paper fills this gap by considering both demand uncertainties. 

The second characteristic of our model compared to existing investigations concerns the multimodality of the demand distribution which is a direct consequence of jointly including probabilistic and fuzzy  uncertainties.
According to \cite{hanasusanto_distributionally_2015}, multimodality of demands is observed, in the following practical situations:  
\begin{itemize}

\item New products: The prediction of the success or failure of new products introduction to the market is difficult. In this case, it is natural to assign the product a bimodal demand distribution.

\item Large customers: In the presence of a large customer that accounts for a large share of the sales, the demand may be multimodal due to irregular bulk orders.

\item	New market entrants: The emergence of a new major competitor can have a significant impact on demand and some customers may switch to this new competitor. 

\item	Fashion trends: These products are typically very volatile and subject to popularity of  some garments. 

\end{itemize}

An important factor that we add to this list is  the impact of social media and customer feedbacks as well as expert judgments on ongoing sales. At the time of inventory replenishment, particularly for newsvendor type products, it may be unclear how to identify the sensitivity of future customers towards experts and past customer  feedbacks. Indeed, being exposed to collective choices is shown to be highly random due to social influences. \cite{ salganik_experimental_2006} reported  an increasing level of unpredictability of market success when customers are exposed to others’ choices.

Another investigation, \cite{rekik_enriching_2017} considered a multimodal form for the demand distribution by assuming two probability functions for the demand subject to a judgment based parameter. The authors assumed that the demand may be one of the two probability distributions weighted with a Bernoulli coefficient. This coefficient, assumed to be deterministic, models the weight given by the decision maker to each probability distribution and  results in a newsvendor problem with a Gaussian Mixture Model (GMM) assumption for the demand.
Our paper extends the model of  \cite{rekik_enriching_2017} by mixing both statistical and fuzzy types of demand data. The subjective uncertainty in the model is modelled with a fuzzy distribution, rather than a deterministic parameter, to better fit with the practical applications described above. 

Furthermore, classical newsvendor models are usually based upon the assumption of risk neutrality (\cite{Khouja1999}; \cite{LeeandNahmias1993}; \cite{Porteus1990}). Recently, there is a growing body of literature that attempts to use alternative risk preferences rather than risk neutrality to describe the newsvendor decision-making behavior (\cite{agrawal_impact_2000},  \cite{guo_newsvendor_2014},  \cite{kamburowski_distribution-free_2014},  \cite{lee_loss-averse_2015}, \cite{wang_loss-averse_2009},  \cite{wang_would_2009}, \cite{Wangetal2012},  \cite{wu_risk-averse_2013}). In this contribution, we also take into account the decision maker attitude towards risk. Indeed, given the introduction of a subjective and human judgment based measure of the demand uncertainty, decision makers decide  replenishment quantities based on their own estimate of the subjective uncertainty, and consequently, they are not neutral to the risks pertaining to their decisions.

Given the prevalence of social influence on demand forecasting and the increasing trend where stock decisions are judgmentally adjusted, we build a hybrid model mixing statistical forecasts and fuzzy weights for a newsvendor problem. Our model contributes to the literature in three levels: 
\begin{itemize}
    \item The demand uncertainty nature by including statistical based and subjective human based sources of demand.
    \item The multi-modality of the demand distribution resulting from mixing probabilistic function and fuzzy weights. 
    \item The risk averse nature of the solution. 
\end{itemize}

The rest of the paper is organized as follows. Section \ref{newsvendorsec} recalls briefly the classical newsvendor formulation. Section \ref{numbers} recalls the basic concepts of  the fuzzy calculus used in this paper. Section \ref{fuzzydistsec} details our newsvendor optimization in a fuzzy Gaussian mixture model permitting to mix probabilistic demand with fuzzy weights. Section \ref{numerical} illustrates the applicability of our model to capture customers' sensitivity to online reviews and provides some numerical insights. Finally, Section \ref{conclusionssec} summarizes the contribution of the paper and presents potential future research directions.

\section{Newsvendor problem}
\label{newsvendorsec}
In this article, we denote the cost per unit for the newsvendor charged by the supplier by $C$,
the cost per unit for the merchant charged by the newsvendor by $M$, and the salvage value per unit by $V$.

If the order quantity is denoted by $Q$, then the profit $\pi$ of 
the newsvendor is a function 
of the demand $X$: 
\begin{equation} \label{newsvendprofit}
\pi(X) =
\begin{cases}
(M-C)X +(V-C)(Q-X),     \quad X \leq Q, \\
(M-C) Q, \quad Q < X,  
\end{cases} 
= 
\begin{cases}
A X +a Q,     \quad X \leq Q, \\
bQ,  \quad Q < X,  
\end{cases} 
\end{equation}
where, for simplicity, we have set 
\begin{equation} \label{AaBb}
A= M-V > 0 , \quad a=V-C < 0 , \quad b= M-C >0. 
\end{equation}

Since the demand $X$ is a random variable, the profit 
$\pi$ is also a random variable, and the newsvendor wishes to adjust 
the order quantity $Q$ in a way that the expectation of the profit is 
maximized. Let us denote the probability distribution function of $X$ by 
$f_X$, and consider its cumulative distribution function by 
$F_X(x) = \int_{-\infty}^x f_X(t) \, dt$. The expectation of the profit is given by 
\begin{eqnarray} \label{profexp}
\bE [ \pi(X) ] &=& 
\int_{-\infty}^Q \big( Ax+aQ \big )  f_X(x) \, dx  + \int_{Q}^\infty 
    bQ    f_X(x) \, dx. 
\end{eqnarray}
For the optimization, by direct calculation one finds that
\begin{equation}
\label{deriviativeprofitexp}
\frac{d}{dQ} \bE [ \pi (X) ] =F_X(Q) (-M+V)+M-C, 
\end{equation}
which is 0, only at the point 
\begin{equation} \label{newsvendorsolution}
Q^*= F_X^{-1} \left ( \frac{M-C}{M-V}\right ). 
\end{equation}
Since  
$
\frac{d^2}{dQ^2} \bE[ \pi(X) ] = 
-f_X(Q) (M -V) < 0,  
$  
the expectation $ \bE [ \pi(X) ] $ is a concave function of $Q$, and it attains its unique maximum at $Q^*$ given by 
\eqref{newsvendorsolution}. 

Here, we derive a formula for the expectation of the profit:
\begin{equation}
\label{profitexpsimpleformula}
\bE [ \pi (X) ] 
= 
(-M+V) \int_0^Q F_X(x) \, dx 
+ (M-C)Q,  
\end{equation}
which can be justified easily by saying that 
the two sides of the equation have the same 
derivative using \eqref{deriviativeprofitexp}, 
and that they both equal to 0 if the order 
quantity $Q$ is 0. Also, the variance of the 
profit can be calculated by the following 
concise formula:
\begin{equation}
\label{profitvarsimpleformula}
\mathrm{Var} [ \pi (X) ] 
= 
\end{equation}
\begin{equation*}
(M-V)^2 \left (
2Q \int_0^Q F_X(x) \, dx - 2 \int_0^Q x F_X(x) \, dx 
- \left ( \int_0^Q F_X(x) \, dx \right)^2 
\right ), 
\end{equation*}
which is used for incorporating risk analysis of 
the newsvendor problem by means of mean-variance analysis, 
see \cite{Choietal2008}, \cite{WuWangetal2009},  
and references therein. That is, their objective is to 
optimize utility functions of the form 
$\bE [ \pi (X) ]  - \beta \, \mathrm{Var} [ \pi (X) ] $ 
and alike, where the decision maker chooses the factor $\beta$ 
depending on their risk attitude. It is easy to see that 
$\frac{d}{dQ} \mathrm{Var} [ \pi (X) ] > 0$, hence exposure 
of the decision maker to more risk by increasing the order quantity. We will see in this article 
that our fuzzy approach will give rise to a novel utility 
function that can naturally incorporate a risk analysis.

\section{Fuzzy numbers and fuzzy random variables}
\label{numbers}

\subsection{Fuzzy numbers, fuzzy expected value, $\al$-cuts, 
and extension principle}

The concept of a {\it fuzzy number} is useful for modelling a 
quantity when there is a range of possibilities for the quantity, 
and one wishes to give values (or weights) 
ranging from 0 to 1 to the possibilities. More precisely, 
a fuzzy number $\widetilde r = (r_1, r_2, r_3, r_4; L, R)$ 
is a function of the form
\[
\widetilde r(x) = 
\begin{cases}
0, \quad x <  r_1,\\
L(x), \quad r_1 \leq x \leq r_2,\\
1, \quad r_2 < x < r_3,\\
R(x), \quad r_3 \leq x \leq r_4,\\
0, \quad r_4 < x.\\
\end{cases}
\]
where $L:[r_1, r_2] \to \R$ is a continuous non-decreasing function, and 
$R: [r_3, r_4] \to \R$ is a continuous non-increasing function such that 
$L(r_1) = R(r_4) = 0$ and $L(r_2) = R(r_3) = 1$. That is, a fuzzy number 
is not a certain number, and allows one to think of $\widetilde r(x)$ to be the 
value of the possibility that the uncertain quantity will be $x$.  
We will denote the set of fuzzy numbers by $\mathcal{F}(\R)$.

If the functions $L$ and $R$ in the definition of a fuzzy number $\widetilde r$
are linear, then $\widetilde r$ is said to be a {\it trapezoidal fuzzy number} and one  
uses the notation $\widetilde r = (r_1, r_2, r_3, r_4)$. Therefore a trapezoidal 
fuzzy number is of the form: 
\[
\widetilde r(x) = 
\begin{cases}
0, \quad x < r_1,\\
\frac{1}{r_2-r_1}(x-r_1), \quad r_1 \leq x \leq r_2,\\
1, \quad r_2 < x < r_3,\\
\frac{1}{r_4-r_3}(r_4-x), \quad r_3 \leq x \leq r_4,\\
0, \quad r_4 < x.\\
\end{cases}
\]

Among a large number of {\it defuzzification} methods in the literature, 
in \cite{LiuLiu2002}, the {\it expected value} of a fuzzy number 
$\widetilde r \in \mathcal{F}(\R)$ is defined as follows. First we 
recall that the {\it possibility}, {\it  necessity} and {\it  credibility} 
of $\{\widetilde r \geq y  \}$ for any $y \in \R$ are defined 
by 
\begin{eqnarray*}
\mathrm{Pos}(\{\widetilde r \geq y  \}) &=& \sup_{x \geq y} \widetilde r(x), \\
\mathrm{Nec}( \{\widetilde r \geq y  \} ) &=& 1 - \mathrm{Pos}(\{\widetilde r < y  \} ) 
= 1 - \sup_{x <y }\widetilde r(x), \\
\mathrm{Cr}( \{ \widetilde r \geq y  \} ) &=& 
\frac{1}{2} \left ( \mathrm{Pos}( \{\widetilde r \geq y  \} ) 
+ \mathrm{Nec}( \{ \widetilde r \geq y  \} ) \right ). 
\end{eqnarray*}
In a similar manner, $ \mathrm{Cr}( \{\widetilde r \leq y  \} )$ is defined. 
Having these notions introduced, in \cite{LiuLiu2002}, the expected 
value of a fuzzy number $\widetilde r$ is defined to be: 
\begin{equation} \label{expgeneral}
E[\widetilde r] = \int_0^\infty \mathrm{Cr}(\{ \widetilde r \geq x \}) \, dx 
-  
\int_{-\infty}^0 \mathrm{Cr}(\{ \widetilde r \leq x \}) \, dx. 
\end{equation}
In Example 2 of \cite{LiuLiu2002}, it is shown that for any trapezoidal fuzzy number 
$\widetilde r=(r_1, r_2, r_3, r_4)$:
\begin{equation} \label{exptrapezoidal}
E[\widetilde r] = E[(r_1, r_2, r_3, r_4)] = \frac{1}{4}(r_1+r_2+r_3+r_4). 
\end{equation}

\begin{figure}
\begin{center}
  \includegraphics[scale=0.4]{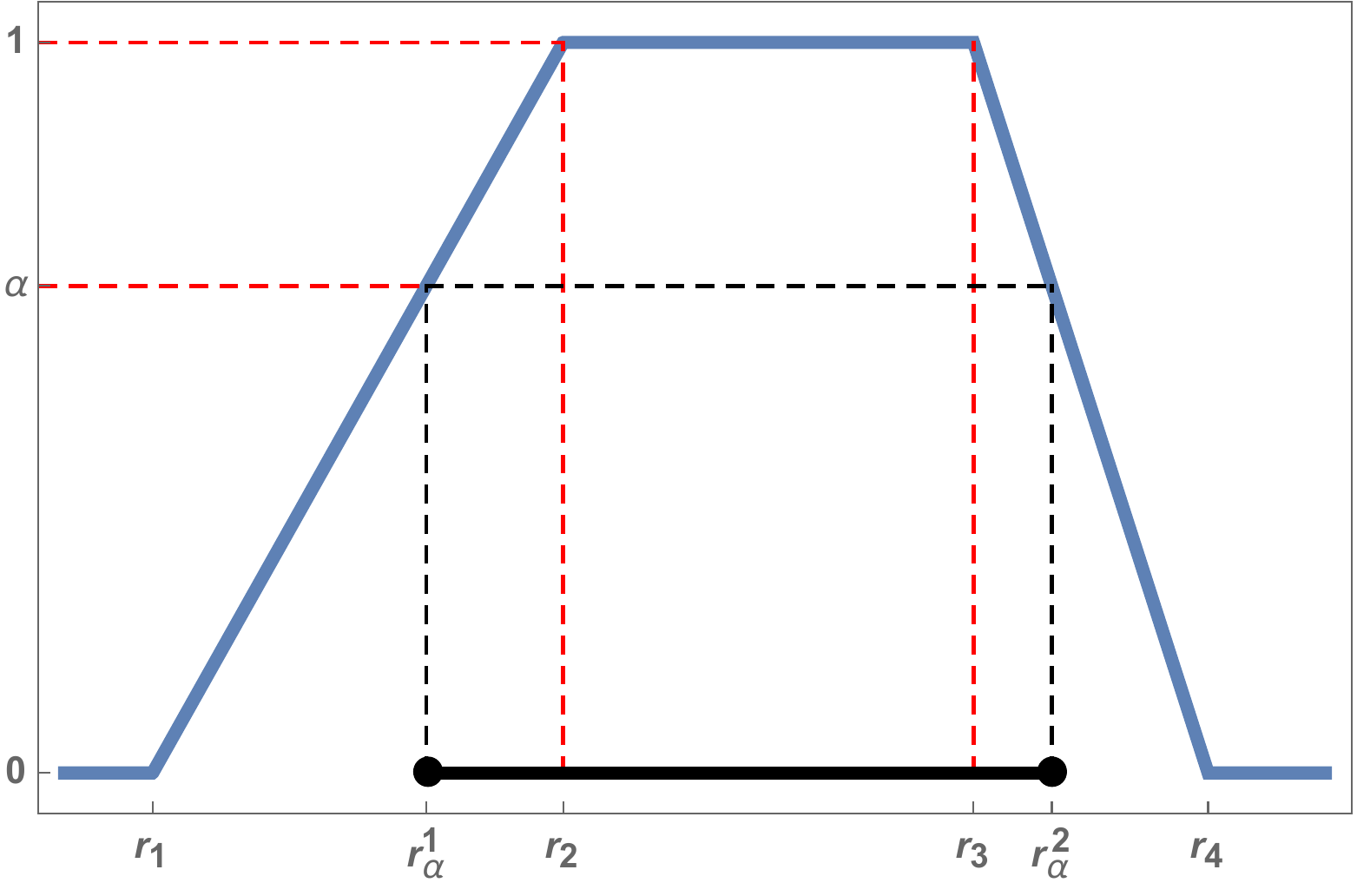}
\end{center}  
  \caption{A trapezoidal fuzzy number $\widetilde r=(r_1, r_2, r_3, r_4)$ 
  and its $\alpha$-cut $[ r^1_\alpha,  r^2_\alpha]$ for some $\alpha$ 
  between 0 and 1.}
  \label{fig:fuzzynum}
\end{figure}

For explicit calculations with fuzzy numbers, It is often handy to use the {\it $\alpha$-cuts} 
of a fuzzy number $\widetilde r \in \mathcal{F}(\R)$. 
By definition, for $0 < \alpha \leq 1$, the $\alpha$-cut of $\widetilde r$ is the closed interval 
$[r^1_\alpha, r^2_\alpha]$ defined by: 
\[
[r^1_\alpha, r^2_\alpha]=\{ x \in \R:   \widetilde r (x) \geq \alpha \}.  
\]
For $\alpha = 0$ the $0$-cut $[r^1_0, r^2_0] $ is defined to be the closed interval $[r_1, r_4]$. 
In Figure \ref{fig:fuzzynum} a schematic representation of a general trapezoidal fuzzy number 
and an $\alpha$-cut is given. In fact, if $\tr = (r_1, r_2, r_3, r_4)$ is a trapezoidal fuzzy 
number, then 
\[
r^1_\alpha = r_1 + (r_2-r_1) \al, \quad r^2_\alpha = r_4 + (r_3-r_4) \al. 
\]

\begin{remark} \label{defuzremark}
For any trapezoidal fuzzy number $\tr = (r_1, r_2, r_3, r_4)$ with the 
$\al$-cuts as above, one can recover the expected value $E[\tr]$ given 
by \eqref{exptrapezoidal} by the following formula:
\[
\frac{1}{2}\left ( \int_0^1 r^1_\al \, d\al + \int_0^1 r^2_\al \, d\al \right ).
\] 
\end{remark}

It is clear that the $\alpha$-cuts as above for $ 0 \leq \alpha \leq 1$ 
determine the fuzzy number $\widetilde r$. Therefore one can 
identify any fuzzy number $\tr \in \mathcal{F}(\R)$  with 
its $\alpha$-cuts and write 
\[
\widetilde r= \{(r^1_\alpha, r^2_\alpha):  0 \leq \alpha \leq 1\}.
\] 
Now, given two fuzzy numbers 
$\widetilde r= \{(r^1_\alpha, r^2_\alpha):  0 \leq \alpha \leq 1\}$ 
and $\widetilde s= \{(s^1_\alpha, s^2_\alpha):  0 \leq \alpha \leq 1\}$ 
in $\mathcal{F}(\R)$, one can define their some $\widetilde r + \widetilde s$ 
and product $\widetilde r \, \widetilde s$ to be the unique fuzzy numbers 
whose representations in terms of 
$\alpha$-cuts are given by 
\[
\widetilde r + \widetilde s 
= 
\{(r^1_\alpha + s^1_\alpha, r^2_\alpha + s^2_\alpha):  0 \leq \alpha \leq 1\}, 
\]
and 
\[
\widetilde r  \, \widetilde s = \{
\left ( \min (r^1_\alpha  s^1_\alpha, r^1_\alpha  s^2_\alpha, r^2_\alpha  s^1_\alpha, r^2_\alpha  s^2_\alpha), 
\max (r^1_\alpha  s^1_\alpha, r^1_\alpha  s^2_\alpha, r^2_\alpha  s^1_\alpha, r^2_\alpha  s^2_\alpha) 
\right):  0 \leq \alpha \leq 1\}. 
\]
We note that 
$\widetilde r \, \widetilde s$ 
$=$  $\{(r^1_\alpha  s^1_\alpha, r^2_\alpha  s^2_\alpha):  0 \leq \alpha \leq 1\}, $
if  $r^1_\alpha, s^1_\alpha \geq 0$. 
Without using the $\alpha$-cuts, the above operations are written in 
terms of the following more complicated formulas: 
\[
( \widetilde r + \widetilde s )(z) = \sup_{x+y =z} \min( \widetilde r (x), \widetilde s(y)), 
\quad 
( \widetilde r \, \widetilde s )(z) = \sup_{x y  =z} \min( \widetilde r (x), \widetilde s(y)). 
\]

Another important tool that we need is the {\it extension principle} 
which allows one to extend the domain of a function defined on real 
numbers to fuzzy numbers. That is, given a function $h: \R \to \R$ and a 
fuzzy number $\widetilde r$ one defines 
$\widetilde s = h(\widetilde r) \in \mathcal{F}(\R)$ by: 
\[
\widetilde s(z) =  h(\widetilde r)(z) = \sup_{h(x)=z} \widetilde r(x). 
\]
In practice, one can conveniently work with the $\alpha$-cuts of 
$\widetilde s = h(\widetilde r) $ which allow one to write: 
\begin{equation} \label{extprinc1}
\widetilde s = h(\widetilde r) = \left \{ ( s^1_\alpha, s^2_\alpha )  : 0 \leq \alpha \leq 1 \right \},
\end{equation}
where
\begin{equation} \label{extprinc2}
s^1_\alpha = \min_{ x \in [r^1_{\alpha}, r^2_{\alpha}] }  h(x),
\qquad 
s^2_\alpha = \max_{ x \in [r^1_{\alpha}, r^2_{\alpha}] }  h(x).    
\end{equation}

\subsection{Fuzzy random variables and their expectations as  fuzzy numbers}

We now discuss the notion of a fuzzy random variable $\widetilde X$. Recall that 
a random variable $X$ on a probability space $(\Omega, \mathcal{M}, P)$, where 
$\mathcal{M}$ is its $\sigma$-algebra of events, is a 
measurable map $X: \Omega \to \R$ and its expectation is defined to be 
$\mathbb{E}[X] = \int_\Omega X(\omega) \, dP(\omega)$.  
The random variable $X$ induces a probability measure $\mu_X$ on $\R$ 
defined by $\mu_X([a, b])$$ = $ $P(X \in [a, b])$$=$$P(X^{-1}([a, b]))$ 
for any interval $[a, b] \subset \R$. If $\mu_X$ is absolutely 
continuous with respect to the Lebesgue measure, one can write 
$d \mu_X(x) = f(x) \, dx$, where $f_X(x)$ is called the probability 
density function of $X$. Therefore $P(X \in [a, b]) = \int_a^b f_X(x) \, dx$, 
and it then follows that $\mathbb{E}[X] = \int_{-\infty}^\infty x f_X(x) \, dx$ 
and more generally $\mathbb{E}[h(X)] = \int_{-\infty}^\infty h(x) f_X(x) \, dx$.

A fundamental fact here is that for any interval $[a, b] \subset \R$, 
and any ordinary random variable $X$ as above, the probability $P(X \in [a, b])$ 
is a certain real number. However, there are situations in which there is no such certainty, 
and one wishes to consider a range of valued possibilities for $P(X \in [a, b])$, where 
the notion of a fuzzy number is useful. That is, one can define a {\it fuzzy random variable} 
$\tX$ on a probability space $(\Omega, \mathcal{M}, P)$ to be a map from $\Omega$ 
to the set of fuzzy numbers $\mathcal{F}(\R)$, denoted in the $\alpha$-cuts notation 
by $\tX(\omega) = \{(X^1_\alpha (\omega),  X^2_\alpha(\omega) ): 0 \leq \alpha \leq 1 \}$, 
such that the functions $X^1_\alpha, X^2_\alpha : \Omega \to \R$ are measurable, and 
the probability of $\tX$ belonging to an interval $[a, b] \subset \R$ is a fuzzy 
number given by $\alpha$-cuts as:
\[
\widetilde P(\tX \in [a, b]) 
= 
\Big \{ \left ( \int_a^b  f^1_{\alpha}(x) \, dx, \int_a^b f^2_\alpha(x) \, dx \right ): 0 \leq \alpha \leq 1  \Big \}, 
\] 
for suitable probability distribution functions $f^1_\alpha$ and $f^2_\alpha$. 
Therefore, the {\it expectation} $\widetilde \bE [ \tX] = \{(\bE[X^1_\al], \bE[X^2_\al]) : \alzeroone\}$ 
is also a fuzzy number, which can be calculated by writing 
\[
\widetilde \bE [ \tX] 
=\Big \{ \left ( \int_{-\infty}^\infty   x f^1_\alpha(x) \, dx, \int_{-\infty}^\infty x f^2_\alpha(x) \, dx \right ): 
0 \leq \alpha \leq 1  \Big \}.
\]

\section{The newsvendor optimization in fuzzy Gaussian mixture model}
\label{fuzzydistsec}

In the modeling of the distribution of the demand in the newsvendor 
problem, considering the technological advancements in recent years, 
there are two main types of data that are important to be considered, namely 
the online reviews and the historical data. Therefore, we wish to 
consider a {\it Gaussian mixture model} for the distribution of the demand of the form 
\[
X \sim p f_1(x) + (1-p) f_2(x), 
\]
where $ 0 \leq p \leq 1$, and for each $i=1,2$, the function 
$f_i(x)$ is the probability density function of the normal distribution with 
mean $\mu_i$ and standard deviation $\sigma_i$, namely:
\begin{equation} \label{normaldist12}
f_i(x)= (2 \pi \sigma_i^2)^{-1/2} 
\exp \left ( - (x- \mu_i )^2/ (2 \sigma_i^2)   \right ), \quad i=1, 2. 
\end{equation} 
We note that we shall denote by $F_i$ the cumulative distribution functions 
corresponding to each $f_i$. 
 The $\mu_i$ and $\sigma_i$  are determined respectively for $i=1,2$ 
by considering the online reviews and the historical data. 
Since one has access to large amounts of each of these data types, it is natural 
to consider them to have normal distributions and to derive values for the $\mu_i$ 
and $\sigma_i$ from the available data. 
However, there is considerable uncertainty 
in choosing the weight $p$ appropriately since one does not know  how the two types 
of information come together to form the forthcoming demand 
$X$.  Therefore in this article we use fuzzy numbers to analyze models that incorporate the 
uncertainties surrounding the weight $p$.

\subsection{Gaussian mixture model with $p = E[\tp]$ for the demand $X$} 
\label{1stmodelsubsec}

In this model, we use the online reviews and the historical data (along with 
information such as the availability of the reviews to potential merchants, their 
social habits, decision making behavior, etc.) to find a trapezoidal fuzzy 
number $\tp = (p_1, p_2, p_3, p_4) $ (with $0 \leq p_1 \leq p_2 \leq p_3 \leq p_4 \leq 1$) 
that incorporates uncertainties regarding how the two data types weigh in the formation of 
the demand $X$. We then consider the probability distribution function of the demand $X$ 
to be of the form 
\begin{equation} \label{probdistpEtp}
f_X(x) = E[\tp] f_1(x) + (1 - E[\tp] ) f_2(x), 
\end{equation}
where $E[\tp] = (p_1 + p_2 + p_3 +p_4)/4$ is the expected value 
of $\tp$ introduced by \eqref{exptrapezoidal}, and,  for each $i=1, 2$, $f_i$ is as in \eqref{normaldist12}. 

In this case, one can use the classical solution of the newsvendor problem as discussed in 
\S \ref{newsvendorsec} to find the optimal order quantity. That is, one can 
use \eqref{newsvendorsolution} to optimize the expectation of the profit 
by ordering the quantity 
\begin{equation} \label{optimalQpEtp}
Q_{\tp}^*= F_X^{-1} \left ( \frac{M-C}{M-V}\right ), 
\end{equation}
where 
the role of the fuzzy number $\tp$ is that its expected value $E[\tp]$ appears in 
$F_X(x) = \int_{-\infty}^x  \left ( E[\tp] f_1(t) + (1 - E[\tp] ) f_2(t) \right ) \, dt$ $=$ 
$   E[\tp] F_1(x) + (1 - E[\tp] ) F_2(x)  $.

\subsection{Gaussian mixture model for $X$ with $p_m \leq p \leq p_M$ and optimization} 
\label{2ndmodelsubsec}

We now  consider the case when one wishes to consider the demand $X$ to have a Gaussian 
mixture distribution where the weigh $p$ is uniformly distributed in a subinterval $[p_m, p_M]$ of 
the interval $[0, 1]$. Therefore, one has  
\begin{equation} \label{densityuniformp}
f_X(x) = p f_1(x) + (1 - p ) f_2(x),
\end{equation}
with $f_i$ for $i=1, 2$ as in \eqref{normaldist12}, and the weight is a {\it rectangular} fuzzy number 
in the sense that it is a special trapezoidal fuzzy number of the form $(p_m, p_m, p_M, p_M)$.

For this model, the expectation of the profit $\bE[\pi(X)]$ 
is clearly a function of the weight $p$. 
Thus, we find the optimal order quantity $Q$ that maximizes 
the average defined by 
\begin{equation} \label{aveprofitexpuniformp}
\overline \bE[\pi(X)] = \frac{1}{p_M-p_m} \int_{p_m}^{p_M} \bE[\pi(X)](p) \, dp. 
\end{equation}
Using  \eqref{profexp} and \eqref{densityuniformp}, we have 
\begin{eqnarray}
 \bE[\pi(X)](p) &=& p \left ( \int_{-\infty}^Q (A x +aQ) f_1(x) \, dx +  \int_{Q}^\infty bQ  f_1(x) \, dx  \right ) \nonumber \\
 &+& (1-p) \left ( \int_{-\infty}^Q(Ax+ a Q) f_2(x) \, dx +  \int_{Q}^\infty b Q f_2(x) \, dx  \right ). \nonumber 
\end{eqnarray}
The latter yields: 

\begin{eqnarray}
\overline \bE[\pi(X)] &=& \frac{1}{p_M- p_m}\int_{p_m}^{p_M} \bE[\pi(X)](p) \, dp \nonumber 
\end{eqnarray}
{\small
\begin{eqnarray}
&=& \frac{1}{2} \left(p_M+p_m\right)
 \Big ( \int_{-\infty}^Q A x  f_1(x) \, dx + a Q F_1(Q) 
 + 
 bQ \left (1-F_1(Q) \right )  \Big) \nonumber \\
 &+& \frac{1}{2}  \left(2-p_m-p_M\right) \Big  ( \int_{-\infty}^QAx  f_2(x) \, dx+ a Q F_2(Q)  
   + b Q \left (1-F_2(Q) \right )  \Big ). \nonumber 
\end{eqnarray}
}

Now we can optimize the average $\overline \bE[\pi]$ by directly calculating 
its first and second derivatives of $Q$. By direct calculation, and replacing 
$A, a, b$ with the values given in \eqref{AaBb}, we find that: 
\[
\frac{d}{dQ} \overline \bE[\pi(X)] = 
\]
\[
\frac{1}{2} \Big (\big (\left(p_m+p_M-2\right) F_2(Q) -
\left(p_m+p_M\right) F_1(Q)\big ) (M-V)-2 (C-M)\Big). 
\]
Therefore, the only solution of the equation $\frac{d}{dQ} \overline \bE[\pi(X)] =0$ 
is the point 
\begin{equation} \label{optimalQuniformp} 
\overline Q_{[p_m, p_M]}^* =  \left ( \left ( \frac{p_m+ p_M}{2}\right ) F_1 + 
\left ( 1- \frac{p_m+ p_M}{2} \right ) F_2 \right )^{-1}  \left ( \frac{M-C}{M-V}\right ). 
\end{equation}
We also find that the second derivative is given by
\[
\frac{d^2}{dQ^2} \overline \bE[\pi(X)] = \frac{1}{2} \left(\left(p_m+p_M\right) \left(f_2(Q)-f_1(Q)\right)-2 f_2(Q)\right) (M-V). 
\]
It can now be seen that $\frac{d^2}{dQ^2} \overline \bE[\pi(X)] < 0$, thus 
$\overline \bE[\pi(X)]$ is a concave function of $Q$. Because, we clearly 
have the natural condition $M-V>0$, and we write the other 
factor $\left(\left(p_m+p_M\right) \left(f_2(Q)-f_1(Q)\right)-2 f_2(Q)\right) $ 
of $\frac{d^2}{dQ^2} \overline \bE[\pi(X)]$ as 
the following to see that it is the sum of three negative terms: 
\[
(p_m -1)f_2(Q) - (p_m+p_M) f_1(Q)+(p_M-1)f_2(Q) <0. 
\]
Therefore $\frac{d^2}{dQ^2} \overline \bE[\pi(X)] < 0$, and the optimal quantity 
$\overline Q_{[p_m, p_M]}^*$ given by \eqref{optimalQuniformp}  uniquely  
 maximizes the average profit expectation $\overline \bE[\pi(X)]$ in this model 
 given by \eqref{aveprofitexpuniformp}.

 We now easily observe a compatibility between the optimal 
 quantities $Q_{\tp}^*$  given by \eqref{optimalQpEtp} for the model 
 \eqref{probdistpEtp}, and the 
 optimal quantity $\overline Q_{[p_m, p_M]}^*$ given by \eqref{optimalQuniformp} 
 for the model  \eqref{densityuniformp}. That is, since $E[(p_m, p_m, p_M, p_M)] 
 = (1/4)(2p_m+2p_M)=(1/2)(p_m+p_M)$, we have: 
 \begin{equation} \label{compatible12}
 Q_{(p_m, p_m, p_M, p_M)}^*= \overline Q_{[p_m, p_M]}^*. 
 \end{equation}

\subsection{Fuzzy Gaussian mixture model for $X$ and fuzzy optimization} 
\label{3rdmodelsubsec}

We now consider a model for the demand to be a fuzzy random variable 
of Gaussian mixture type and will denote it $\tX$. The significance of this 
model is that it can incorporate a general trapezoidal fuzzy number $\tp$ 
as the weight of the online reviews that contribute to the formation of the 
demand. We then use the machinery of fuzzy calculus to calculate the 
the profit $\tY = \pi(\tX)$ as a fuzzy random variable whose expectation 
$\bE(\tY)$ will be a fuzzy number with non-linear sides: it will be more 
general than the trapezoidal case as the end points of its $\al$-cuts 
are given by quadratic expressions in $\al$. We shall then defuzzify 
$\bE(\tY)$ in a natural manner for our management purposes, and 
will present the optimal order quantity $Q$ that maximizes the defuzzified 
quantity.  

Therefore, we consider the normal probability distributions $f_1$ and $f_2$ 
as in \eqref{normaldist12}, which separately reflect the distribution of the 
demand based on the online reviews and historical sales.  We then use the 
information such as the expert opinion, availability of the reviews to potential 
merchants,  social habits of the merchants, decision making behavior, etc., 
to find a fuzzy number 
\[
\tp =(p_1, p_2, p_3, p_4), \quad  0 \leq p_1 \leq p_2 \leq p_3 \leq p_4 \leq 1, 
\]
which reflects the weight of the online reviews in the formation of the demand 
and can incorporate the uncertainties in the weight. In our model, we choose the 
numbers $p_1$ and $p_2$  independently from the numbers $p_3$ and $p_4$ by using 
separate parts of the data so that the random variables $X^{1}_\alpha$ and $X^{2}_\alpha$ 
whose probability distributions are defined below are independent: 
\[
f_{X^{1}_\alpha}(x) = p_{1}^\alpha f_1(x) + (1- p_{1}^\alpha) f_2(x),    
\] 
\[
f_{X^{2}_\alpha}(x) = p_{2}^\alpha f_1(x) + (1- p_{2}^\alpha) f_2(x).  
\] 
We now consider the demand to have the fuzzy random variable $\tX$ 
defined in terms of its $\al$-cuts by
\begin{equation} \label{fuzzydemandintro}
\tX = \{ (X^{(1)}_\al, X^{(2)}_\al), \alzeroone \},
\end{equation}
where 
\[
X^{(1)}_\alpha = \min (X^{1}_\alpha, X^{2}_\alpha), \quad X^{(2)}_\alpha= \max (X^{1}_\alpha, X^{2}_\alpha). 
\]

\begin{lemma} \label{fX12alpha}
The joint probability distribution function of the random variables $X^{(1)}_\alpha$ and 
$X^{(2)}_\alpha$, denoted by $f_{X^{(1,2)}_\al}$, is given by 
\[
f_{X^{(1,2)}_\al}(x_1, x_2) = 
\begin{cases}
f_{X^{1}_\alpha}(x_1) f_{X^{2}_\alpha}(x_2) + 
f_{X^{2}_\alpha}(x_1) f_{X^{1}_\alpha}(x_2), \quad x_1 \leq x_2,  \\
0, \quad x_2 < x_1.
\end{cases}
\] 
\begin{proof}
From the definition of  $X^{(1)}_\alpha$ and $X^{(2)}_\alpha $, we always have 
$X^{(1)}_\alpha \leq X^{(2)}_\alpha $ and therefore $f_{X^{(1,2)}_\al}(x_1, x_2)=0$ when 
$x_2 < x_1$. In order to calculate $f_{X^{(1,2)}_\al}(x_1, x_2)$ 
when $x_1 \leq x_2$, we calculate  
\begin{eqnarray*}
F_{X^{(1,2)}_\al}(x_1, x_2) &=& P\left (X^{(1)}_\alpha \leq x_1, X^{(2)}_\alpha \leq x_2 \right ) \\
&=& P \left (  X^{1}_\alpha \leq x_1, X^{2}_\alpha \leq x_1 \right )
+ P \left (  X^{1}_\alpha \leq x_1, x_1 < X^{2}_\alpha \leq x_2 \right ) \\
&+&P \left (  X^{2}_\alpha \leq x_1, x_1 < X^{1}_\alpha \leq x_2 \right )\\
&=&
F_{X^{1}_\alpha}(x_1) F_{X^{2}_\alpha}(x_1) 
+  F_{X^{1}_\alpha}(x_1) \left (  F_{X^{2}_\alpha}(x_2) - F_{X^{2}_\alpha}(x_1)  \right ) 
\\
&+&  F_{X^{2}_\alpha}(x_1) \left (  F_{X^{1}_\alpha}(x_2) - F_{X^{1}_\alpha}(x_1)  \right ), 
\end{eqnarray*}
which yields: 
\[
f_{X^{(1,2)}_\al}(x_1, x_2) 
= 
\frac{\partial^2 F_{X^{(1,2)}_\al}}{\partial x_1 \partial x_2} (x_1, x_2)
= f_{X^{1}_\alpha}(x_1) f_{X^{2}_\alpha}(x_2) + f_{X^{2}_\alpha}(x_1) f_{X^{1}_\alpha}(x_2). 
\]

\end{proof}
\end{lemma}

We need the following statement for defuzzifications in the sequel, which involve integrations 
over $\al$ from 0 to 1.

\begin{lemma} \label{fX12alphaopenlemma}
When $x_1 \leq x_2$, we have 
\[
\int_0^1 f_{X^{(1,2)}_\al}(x_1, x_2) \, d\alpha = 
\]
\[
P_1 \, f_1\left(x_1\right) f_1\left(x_2\right)+P_2 \left ( f_1\left(x_1\right) f_2\left(x_2\right)+ 
f_2\left(x_1\right) f_1\left(x_2\right) \right)+P_3 \, f_2\left(x_1\right) f_2\left(x_2\right), 
\]
where 
\begin{eqnarray} \label{P1P2P3}
P_1 &=&\frac{p_1 p_3}{3}+\frac{2 p_2 p_3}{3}+\frac{2 p_1 p_4}{3}+\frac{p_2 p_4}{3},  \\
P_2&=&\frac{p_1}{2}-\frac{p_3 p_1}{3}-\frac{2 p_4 p_1}{3}+\frac{p_2}{2}+\frac{p_3}{2}+
\frac{p_4}{2}-\frac{2 p_2 p_3}{3}-\frac{p_2 p_4}{3}, \nonumber \\
P_3 &=& \frac{p_3 p_1}{3}+\frac{2 p_4 p_1}{3}-p_1-p_2+\frac{2 p_2 p_3}{3}-p_3+
\frac{p_2 p_4}{3}-p_4+2. \nonumber
\end{eqnarray}
\begin{proof} 
Using Lemma \ref{fX12alpha} and direct substitutions, when $x_1 \leq x_2$, 
we have: 
\begin{equation*} \label{fX12alphaopen}
f_{X^{(1,2)}_\al}(x_1, x_2) = 
\end{equation*} 
\begin{center}
\begin{math}
2 f_1\left(x_1\right) f_1\left(x_2\right) p^1_{\alpha} p^2_{\alpha}+f_2\left(x_1\right) f_1\left(x_2\right) p^1_{\alpha}-2 f_2\left(x_1\right) f_1\left(x_2\right) p^1_{\alpha} p^2_{\alpha}+f_2\left(x_1\right) f_1\left(x_2\right) p^2_{\alpha}+f_1\left(x_1\right) f_2\left(x_2\right) p^1_{\alpha}-2 f_1\left(x_1\right) f_2\left(x_2\right) p^1_{\alpha} p^2_{\alpha}+f_1\left(x_1\right) f_2\left(x_2\right) p^2_{\alpha}-2 f_2\left(x_1\right) f_2\left(x_2\right) p^1_{\alpha}+2 f_2\left(x_1\right) f_2\left(x_2\right) p^1_{\alpha} p^2_{\alpha}-2 f_2\left(x_1\right) f_2\left(x_2\right) p^2_{\alpha}+2 f_2\left(x_1\right) f_2\left(x_2\right).
\end{math}
\end{center}
Then, we find the result by the substitution 
\[
p^1_\alpha = p_1 + (p_2-p_1) \al, \quad p^2_\alpha = p_4 + (p_3-p_4) \al, 
\]
and direct integration over $\al$ from 0 to 1. 
\end{proof}
\end{lemma}

Also, for the optimizations in the following subsections, we shall need the 
following properties of $P_1, P_2, P_3$ given by 
\eqref{P1P2P3}. 
\begin{lemma} \label{P1P2P3lem}
We have: 
\[
0 \leq P_1,  P_3 \leq 2, \quad 0 \leq P_2 \leq 1, \quad \frac{1}{2}P_1+ P_2+ \frac{1}{2}P_3 = 1.
\]
and 
\[
P_1+2 P_2+ P_3 = 2. 
\]
\begin{proof}
Since $0 \leq p_1 \leq p_2 \leq p_3 \leq p_4 \leq 1$, it is clear that $0 \leq P_1 \leq 2$. 
One can then see that by setting $q_i = 1 - p_i$ for $i=1, \dots, 4$ one has 
\[
P_3 = \frac{q_1 q_3}{3}+\frac{2 q_2 q_3}{3}+\frac{2 q_1 q_4}{3}+\frac{q_2 q_4}{3}, 
\]
which shows that $0\leq P_3 \leq 2$. The identity $P_1 + 2 P_2 +P_3 =2$ can 
be checked directly, which yields: 
\[
P_2 = 1 - \frac{1}{2}P_1- \frac{1}{2}P_3 \leq 1. 
\]
In order to see that $P_2 \geq 0$, we write 
\begin{eqnarray*}
2 P_2 &=& p_1+p_2+p_3+p_4-\frac{2}{3}  p_3 p_1-
\frac{4}{3} p_4 p_1 -\frac{4}{3} p_2 p_3 -\frac{2}{3} p_2 p_4 \\
&\geq& 
p_1^2+p_2^2+p_3^2+p_4^2-\frac{2}{3}  p_3 p_1 -\frac{4}{3} p_4 p_1 -
\frac{4}{3} p_2 p_3 -\frac{2}{3} p_2 p_4 \\
&=&  \begin{pmatrix}
p_1 & p_2 & p_3 & p_4\\
\end{pmatrix} 
\left(
\begin{array}{cccc}
 1 & 0 & -1/3 & -2/3 \\
 0 & 1 & -2/3 & -1/3\\
 -1/3 & -2/3 & 1 & 0 \\
 -2/3 & -1/3 & 0 & 1 \\
\end{array}
\right)
\begin{pmatrix}
p_1 \\
p_2 \\
p_3 \\
p_4
\end{pmatrix}.    
\end{eqnarray*}
The latter is non-negative since the above square matrix is 
symmetric and its eigenvalues are the non-negative numbers 
$2, 4/3, 2/3, 0$. 

\end{proof}
\end{lemma}

By using \eqref{extprinc1} and \eqref{extprinc2}, we apply the 
profit function $\pi$ to the fuzzy random variable $\tX$ and obtain 
the profit of the newsvendor as a fuzzy random variable $\tY = \pi (\tX)$: 
\[
\tX = \left \{ (X^{(1)}_\alpha, X^{(2)}_\alpha):  \alzeroone \right \}, 
\quad \tY = \pi (\tX) = \left \{ (Y^1_\alpha, Y^2_\alpha):  \alzeroone \right \}, 
\]
where
\[
Y^1_\alpha = h_1(X^{(1)}_\al, X^{(2)}_\al) = 
\begin{cases}
 A X^{(1)}_\al + a Q         , \quad  X^{(1)}_\al \leq X^{(2)}_\al \leq Q,\\
    bQ       , \quad Q \leq X^{(1)}_\al \leq X^{(2)}_\al, \\
 \min \Big( A X^{(1)}_\al + a Q,    bQ \Big ),  \quad X^{(1)}_\al \leq Q \leq  X^{(2)}_\al,           
\end{cases}
\]
\[
Y^2_\alpha = h_2(X^{(1)}_\al, X^{(2)}_\al) = 
\begin{cases}
 A X^{(2)}_\al + a Q            , \quad  X^{(1)}_\al \leq X^{(2)}_\al \leq Q,\\
 bQ         , \quad Q \leq X^{(1)}_\al \leq X^{(2)}_\al, \\
(A+a)Q= bQ         , \quad X^{(1)}_\al \leq Q \leq  X^{(2)}_\al.            
\end{cases}
\]

\subsubsection{Optimization of $\int_0^1 \bE[Y^1_\al] \, d\al$}

We have: 
\begin{eqnarray} \label{preexpy1S0}
\bE[Y^1_\al] 
&=& \int_{-\infty}^Q dx_1 \int_{x_1}^Q dx_2\,  \left ( Ax_1 + a Q \right ) f_{X^{(1,2)}_\al}(x_1, x_2)  \\
& +& \int_Q^\infty dx_1 \int_{x_1}^\infty dx_2  \,b Q\,  f_{X^{(1,2)}_\al}(x_1, x_2) \nonumber  \\
 &+& \int_{-\infty}^Q dx_1 \int_Q^{\infty } dx_2 \, \left ( A x_1 + a Q \right ) f_{X^{(1,2)}_\al}(x_1, x_2).  \nonumber
 \end{eqnarray}

In order to find the quantity $Q$ that maximizes $\int_0^1 \bE[Y^1_\al] \, d\al $, 
we use Lemma \ref{fX12alphaopenlemma} and calculate its first and second 
derivatives with respect to $Q$. We find that: 

\begin{eqnarray} \label{firstderleft}
\frac{d}{dQ} \int_0^1 \bE[Y^1_\al] \, d\al  &=& 
 -C+M + \left (-2C + M+V \right)  G(Q)\\
 &+&\left (C-V \right)   H(Q) + \left(-M+V \right) J(Q), \nonumber
 \end{eqnarray}
 where
 \begin{eqnarray} \label{GgHJ}
 G(x)&=&\int_{-\infty}^x g (t) \, dt, \\
 g(t) &=& P_1\, f_1(t) F_1(t) + P_2 \big (f_1(t) F_2(t)+f_2(t) F_1(t) \big ) + P_3 \, f_2(t) F_2(t), \nonumber \\
 H(x) &=&  P_1 \, F_1(x)^2 + 2 P_2 \, F_1(x) F_2(x) + P_3 \, F_2(x)^2,  \nonumber \\
 J(x) &=& P_1 \,F_1(x)+P_2 \left ( F_1(x) + F_2(x)\right )+P_3 \,F_2(x).  \nonumber
\end{eqnarray}

We now make some observations about the functions $G, H, J$, which 
shall be insightful in the optimizations in the sequel.  
\begin{lemma} \label{GHJprop}
The functions $G, H, J$ are non-negative, non-decreasing, for any real number $x$: 
\begin{equation} \label{GhalfH}
G(x) = \frac{1}{2}H(x), \quad H(x) \leq J(x),
\end{equation}
and
\[
\lim_{x \to - \infty} G(x) = \lim_{x \to - \infty} \frac{1}{2}H(x)  =\lim_{x \to - \infty} \frac{1}{2} J(x) = 0, 
\]
\[
\lim_{x \to  \infty} G(x) = \lim_{x \to  \infty} \frac{1}{2} H(x)  =\lim_{x \to  \infty} \frac{1}{2} J(x) = 1. 
\]
\begin{proof}
It is clear that these functions are non-negative and non-decreasing. 
Using the definitions of the functions, we can write
\begin{eqnarray*}
G(x)&=&\int_{-\infty}^x \left (  P_1\, f_1(t) F_1(t) + P_2 \big (f_1(t) F_2(t)+f_2(t) F_1(t) \big ) + P_3 \, f_2(t) F_2(t) \right ) dt\\
&=& \int_{-\infty}^x \frac{d}{dt} \left (   \frac{P_1}{2} F_1(t)^2 + P_2  \,F_1(t) F_2(t)  + \frac{P_3}{2} F_2(t)^2 \right ) dt\\
&=& \frac{P_1}{2} F_1(x)^2 + P_2  \,F_1(x) F_2(x)  + \frac{P_3}{2} F_2(x)^2 \\
&=& \frac{1}{2}H(x). 
\end{eqnarray*}

The fact that 
$H(x) \leq J(x)$ is an easy consequence of the fact that $ 0 \leq F_1(x) \leq 1,$ $0 \leq F_2(x) \leq 1$ 
and that $ 0 \leq P_1, P_2, P_3$ (proved in Lemma \ref{P1P2P3lem}). The vanishing of 
these functions in limit at $-\infty$ is also easily seen, and the limits at $\infty$ are 
calculated by considering the fact that  $P_1/2 + P_2 + P_3/2 =1$.  
\end{proof}
\end{lemma}

We also have  
\[
\lim_{Q \to - \infty} \left ( \frac{d}{dQ} \int_0^1 \bE[Y^1_\al] \, d\al \right )= -C+M >0, 
\] 
\[
\lim_{Q \to  \infty} \left (  \frac{d}{dQ} \int_0^1 \bE[Y^1_\al] \, d\al  \right )=  - C +V < 0. 
\]

For the second derivative we find that:  
\begin{eqnarray} \label{secondderleft}
&&\frac{d^2}{dQ^2} \int_0^1 \bE[Y^1_\al] \, d\al \\ 
&=& 
(V-M) 
\Big (
f_1(Q) \left ( (P_1 (1 - F_1(Q) )+P_2 (1-F_2(Q)) \right) \Big ) \nonumber  \\ 
&+& (V-M) \Big ( f_2(Q) \left(P_2 (1-F_1(Q) ) +P_3 ( 1- F_2(Q) )\right )
\Big) < 0, \nonumber 
\end{eqnarray}
which shows that we have a concave function of $Q$. Therefore, the 
unique optimal quantity $Q$ that maximizes $\int_0^1 \bE[Y^1_\al] \, d\al $ 
is the solution of setting \eqref{firstderleft} equal to 0. 
By defining the function 
\begin{equation} \label{leftcommfunctionS0}
F_{\tp, L}(x) =  \frac{2C - M-V}{M-V} G(x)+ \frac{-C+V}{M-V}  H(x)  + J(x) = -\frac{1}{2}H(x)+J(x),  
\end{equation}
where the last identity follows from  \eqref{GhalfH}, 
this quantity is given by 
\begin{equation}  \label{optimalQleftlegS0}
Q^*_{\tp, L} = F^{-1}_{\tp, L} \left ( \frac{M-C}{M-V}\right ). 
\end{equation}

\subsubsection{Optimization of $\int_0^1 \bE[Y^2_\al] \, d\al$} Similarly, 
we start by writing the following expression: 
 
\begin{eqnarray} \label{preexpy2}
 \bE[Y^2_\al] 
 &=& 
 \int_{-\infty}^\infty dx_1 \int_{x_1}^\infty dx_2 \, 
 h_2(x_1, x_2) \, f_{X^{(1,2)}_\al}(x_1, x_2) 
  \end{eqnarray}
 \begin{eqnarray*}
&=& \int_{-\infty}^Q dx_1 \int_{x_1}^Q dx_2\,  \left ( A x_2 + a Q \right ) f_{X^{(1,2)}_\al}(x_1, x_2) \nonumber \\
&+&  \int_Q^\infty dx_1 \int_{x_1}^\infty dx_2  bQ   f_{X^{(1,2)}_\al}(x_1, x_2) \nonumber \\
&+&   \int_{-\infty}^Q dx_1 \int_{Q}^\infty dx_2   b Q \,  f_{X^{(1,2)}_\al}(x_1, x_2). 
 \end{eqnarray*}

We then use Lemma \ref{fX12alphaopenlemma},  and calculate the 
derivatives of order one and two of $\int_0^1 \bE[Y^2_\al] \, d\al$ 
with respect to $Q$. We find that:   

\begin{eqnarray} \label{firstderright}
\frac{d}{dQ} \int_0^1 \bE[Y^2_\al] \, d\al  &=& -C+M +(-2 C+M+V) G(Q) \\
&+& (C-M) H(Q)  \nonumber 
\end{eqnarray}
where
the functions $G, H, J$ are already given in \eqref{GgHJ}. 
We also note that 
\[
\lim_{Q \to -\infty} \left ( \frac{d}{dQ} \int_0^1 \bE[Y^2_\al] \, d\al \right ) = -C+M > 0, 
\]
\[
\lim_{Q \to \infty} \left ( \frac{d}{dQ} \int_0^1 \bE[Y^2_\al] \, d\al \right ) = -C+V < 0. 
\]

The function $\int_0^1 \bE[Y^2_\al] \, d\al$ is concave in $Q$
 because its second derivative with respect to $Q$ 
is the sum of three negative terms: 
\begin{eqnarray} \label{secondderright}
&&\frac{d^2}{dQ^2} \int_0^1 \bE[Y^2_\al] \, d\al \\ 
&=& 
P_1  f_1(Q) F_1(Q) (-M+V)  \nonumber \\
&+&P_2 \Big ( f_2(Q) \big(F_1(Q) (-M+V)\big)+
f_1(Q) \big(F_2(Q) (-M+V)\big) \Big ) \nonumber \\
&+&P_3  f_2(Q) F_2(Q) (-M+V)    <0. \nonumber
\end{eqnarray}

Therefore, the unique optimal solution that maximizes $\int_0^1 \bE[Y^2_\al] \, d\al$ 
is the quantity $Q$ that is the solution of setting \eqref{firstderright} equal to 0. This is 
given by 
\begin{equation} \label{optimalQrightlegS0}
Q^*_{\tp, R} = F^{-1}_{\tp, R} \left ( \frac{M-C}{M-V} \right ),   
\end{equation}
where considering \eqref{firstderright}  and \eqref{GhalfH}: 
\begin{equation} \label{rightcommfunctionS0}
F_{\tp, R}(x) = \frac{2 C-M-V}{M-V} \,G(x) + \frac{M-C}{M-V}\, H(x) = \frac{1}{2}H(x).  
\end{equation}

\subsubsection{Optimization of 
$
(1 - \beta)  \int_0^1 \bE[Y^1_\al] \, d\al+ \beta \int_0^1 \bE[Y^2_\al] \, d\al 
$} 
\label{betadefuzoptimizesubsec}
In the previous subsections, we have seen that when the demand 
$\tX$ is a fuzzy random variable as introduced in 
\eqref{fuzzydemandintro}, the expectation of the profit $\tY = \pi (\tX)$ 
is a fuzzy number whose $\al$-cuts are 
given by quadratic expressions in $\al$.  By performing an integration 
over $\al$ and finding the optimal quantity $Q^*_{\tp, L} $
given by \eqref{optimalQleftlegS0}, we have introduced an optimal 
quantity that gives the largest possible average 
for the small possible values for the profit. 
Therefore, the quantity $Q^*_{\tp, L}$ can be used when 
the newsvendor wishes to go with minimum risk by choosing the order 
quantity that gives the security that the average of the small possible values 
of the fuzzy profit expectation is maximized. 
In the other direction, 
the optimal quantity $Q^*_{\tp, R}$ given by \eqref{optimalQrightlegS0} 
maximizes the average over the greatest possible values for the fuzzy 
profit expectation. Therefore, the quantity $Q^*_{\tp, R}$ can be 
used if the newsvendor wishes to create maxim possibility for the 
profit at the price of being in the riskiest possible situation.

Thus, it is natural at this stage to use a real number $\beta$ 
in the interval $[0, 1]$  as a {\it risk factor} so that the newsvendor 
can choose an optimal order quantity depending on the risk that 
they can afford.  That is, we wish to find an optimal quantity 
$Q^*_{\tp, \beta}$ that maximizes 
\begin{equation} \label{betadefuzprofitexp}
(1 - \beta)  \int_0^1 \bE[Y^1_\al] \, d\al+ \beta \int_0^1 \bE[Y^2_\al] \, d\al.  
\end{equation}
Note that the latter defines a defuzzification that generalizes 
the one discussed in Remark \ref{defuzremark}.

Considering the calculations performed for the results in the 
previous subsections, we can find this optimal quantity quite 
easily. The reason is that, using \eqref{firstderleft} and 
\eqref{firstderright}, we have an explicit expression for 
the derivative of  \eqref{betadefuzprofitexp} with respect to $Q$. 
Since, using the information in \eqref{secondderleft} and 
\eqref{secondderright}, we have a concave expression, 
the unique optimal quantity that maximizes the defuzzified 
expectation of the fuzzy profit, given by \eqref{betadefuzprofitexp}, 
is the quantity 
\begin{equation} \label{optimalfuzzybeta}
Q^*_{\tp, \beta} =  F^{-1}_{\tp, \beta} \left ( \frac{M-C}{M-V} \right), 
\end{equation}
where 
\begin{eqnarray} \label{defuzcumulfunctionbeta} 
F_{\tp, \beta}(x) &=& (1-\beta) F_{\tp, L}(x) +  \beta F_{\tp, R}(x) \\
&=& \frac{2\beta -1}{2} H(x) + (1- \beta) J(x) \nonumber \\
&=& \frac{1}{2} H(x) + (1 - \beta ) \left ( J(x) - H(x) \right). \nonumber
\end{eqnarray}
Clearly, the optimal solutions for $\beta = 0$ and $\beta = 1$ match 
with the solutions given by \eqref{optimalQleftlegS0} and \eqref{optimalQrightlegS0}, 
respectively.

The resemblance between the formulas for the optimal solution  
 \eqref{optimalfuzzybeta} 
and the classical newsvendor solution \eqref{newsvendorsolution}
raises the question whether the function 
$F_{\tp, \beta}$ given by 
\eqref{defuzcumulfunctionbeta} is the cumulative distribution 
functions of a random variable. The following 
proposition in fact shows that one can take this point of view. 
\begin{proposition} \label{Fscumulativeprop}
The function $F_{\tp, \beta}$, for any $\beta$ in the interval 
$[0, 1]$, satisfies the properties 
of the cumulative distribution function of a random variable: it is  
non-negative, non-decreasing, tends to 0 
at $-\infty$, and tends to 1 at $\infty$. 
\begin{proof}
This is a direct consequence of  Lemma \ref{GHJprop}, 
\eqref{secondderleft}, \eqref{secondderright}, and the 
fact that $F_{\tp, \beta}(x) =$ $(1-\beta) F_{\tp, L}(x) +  \beta F_{\tp, R}(x)$.  
\end{proof}
\end{proposition}

We now make an observation, which is in complete agreement 
with intuition:

\begin{proposition}  \label{Qincreasinginbeta}
The optimal order quantity $Q^*_{\tp, \beta}$ given by 
\eqref{optimalfuzzybeta} is a non-decreasing function of the 
risk factor $\beta$. 
\begin{proof}
Using the fact $H(x) \leq J(x)$ from Lemma \ref{GHJprop}, 
it follows from the formula 
\[
F_{\tp, \beta}(x) =  \frac{1}{2} H(x) + (1 - \beta ) \left ( J(x) - H(x) \right)
\]
that for any real number $x$, the value $F_{\tp, \beta}(x)$ defines a 
non-increasing function of the risk factor $\beta$. Therefore,  
$Q^*_{\tp, \beta} = F^{-1}_{\tp, \beta} \left ( \frac{M-C}{M-V} \right)$ 
is non-decreasing in $\beta$. 
\end{proof}
\end{proposition}

\begin{figure}
\begin{center}
  \includegraphics[scale=0.5]{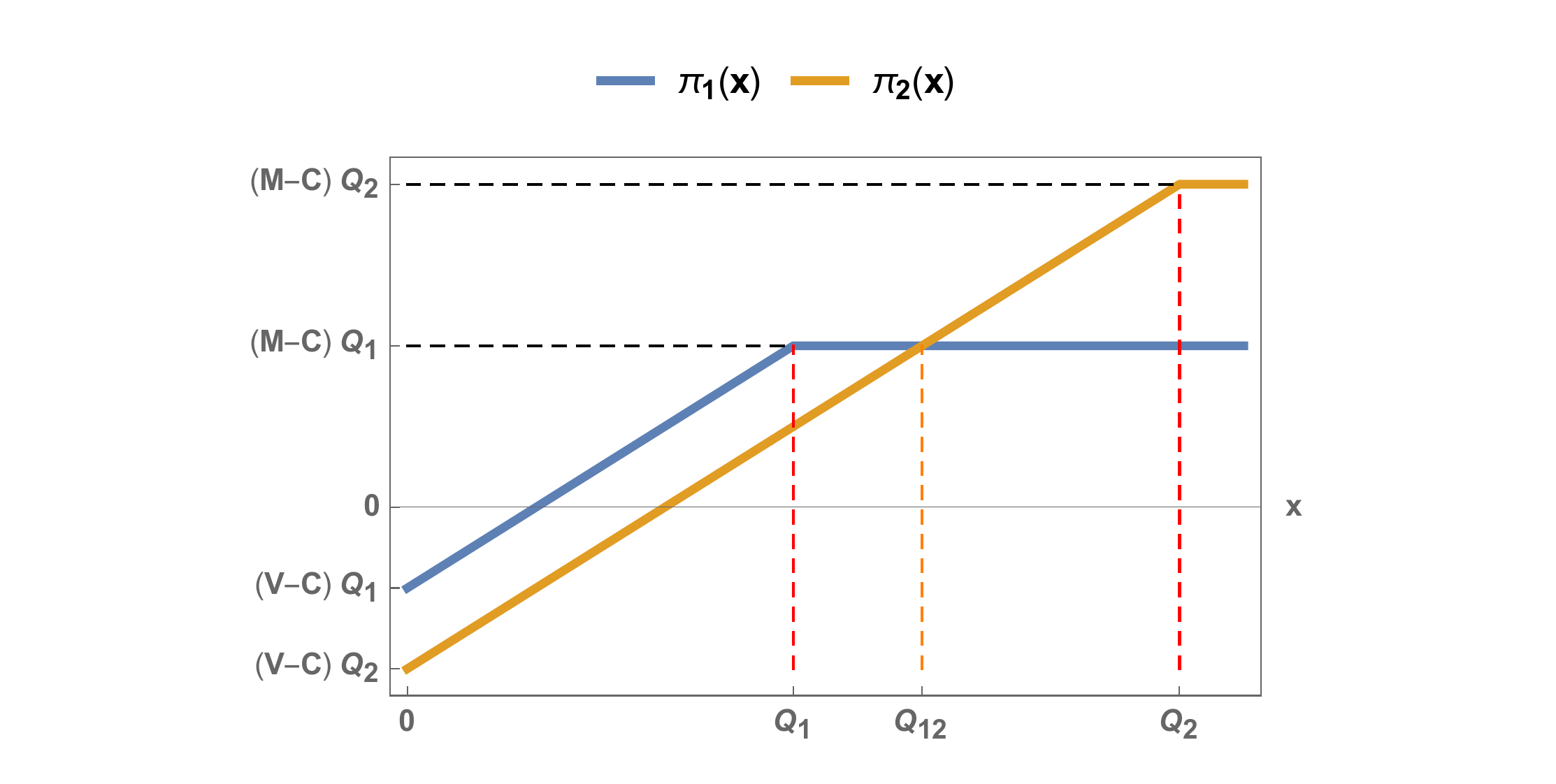}
\end{center}  
  \caption{The profit function $\pi_i$ corresponds to the order quantity $Q_i$ 
  for $i=1, 2$.  }
   \label{fig:Q1Q2}
\end{figure}

The fact that $Q^*_{\tp, \beta}$ is non-decreasing in $\beta$ justifies the reason 
for calling the latter the risk factor. Because when the newsvendor thinks of two 
order quantities $Q_1$ and $Q_2$ such that $Q_1 < Q_2$, the 
profit function $\pi_2$ corresponding to $Q_2$ provides a wider range 
of possibilities for the gain and loss compared to the profit function $\pi_1$ 
corresponding to  $Q_1$, as it can be seen in Figure \ref{fig:Q1Q2}. 
That is,  if the newsvendor chooses the more risky case of ordering the larger 
quantity $Q_2$, then  $\pi_2(x) > \pi_1(x)$ if and only if the demand $x$ is 
greater than  
\[
Q_{12} = \frac{M-C}{M-V} Q_1+\frac{C-V}{M-V}Q_2. 
\] 
Thus, the calculation of the probability that the demand is greater than $Q_{12}$ can help the newsvendor to choose between $Q_1$ and $Q_2$.

\subsection{Newsvendoric defuzzifications of the fuzzy demand} 

Proposition \ref{Fscumulativeprop} confirms that the function  
$F_{\tp, \beta}$ given by  
\eqref{defuzcumulfunctionbeta} defines the cumulative distribution 
function of a random variable, which will be denoted 
by $X_{\tp, \beta}$ hereafter. Hence, given the fuzzy 
random variable $\tX$ defined in \eqref{fuzzydemandintro}, considering 
\eqref{newsvendorsolution}, the optimal quantity given 
by \eqref{optimalfuzzybeta} is  
the solution of the classical newsvendor problem when the demand  
is assumed to be the ordinary random variable  
$X_{\tp, \beta}$, whose distribution density function is:

\begin{eqnarray}  \label{densitytpbeta}
 f_{X_{\tp, \beta}}(x)  &=& \frac{d}{dx} F_{\tp, \beta}(x)   \\
&=& \frac{d}{dx} \left ( \frac{2\beta -1}{2} H(x) + (1- \beta) J(x) \right ) \nonumber 
\end{eqnarray}
\begin{eqnarray*}
&=&
(2 \beta -1) \big (f_1(x) \left(P_1 F_1(x)+P_2 F_2(x)\right)+f_2(x) \left(P_2 F_1(x) +
P_3 F_2(x)\right)\big )   \\
&+&(\beta -1) \left(-\left(P_1+P_2\right) f_1(x)-\left(P_2+P_3\right) f_2(x)\right). 
\end{eqnarray*}

This means that, one can replace the demand $\tX$, which is a 
fuzzy random variable, with an ordinary random variable, 
and consider the solution of the classical newsvendor 
problem to obtain the optimal quantity  
\eqref{optimalfuzzybeta}.  Hence, we have arrived at a natural defuzzification  
of the fuzzy random variable $\tX$ via our fuzzy 
newsvendor optimization.

We made the observation \eqref{compatible12} to show that there is a relation 
between the optimizing quantities for   the models written in \S \ref{1stmodelsubsec} 
and \S \ref{2ndmodelsubsec}. We now make observations in the following remark 
to show that the optimized quantity \eqref{optimalfuzzybeta} and the 
distribution function given by \eqref{densitytpbeta} of the defuzzified 
random variable $\tX_{\tp, \beta}$   relate to the optimized quantities 
\eqref{optimalQpEtp} and \eqref{optimalQuniformp} when $\beta = 1/2$.

First, we explain how the intuitive model of \S \ref{1stmodelsubsec}  can be 
viewed as a particular case of the fuzzy model. 

\begin{remark} \label{model31relate}
Starting with a trapezoidal fuzzy number $\tp = (p_1, p_2, p_3, p_4)$ and the 
fuzzy demand $\tX$ as in \eqref{fuzzydemandintro}, if one sets the risk 
factor $\beta$ equal to $1/2$, then  the density of the 
corresponding defuzzified demand $X_{\tp, \frac{1}{2}}$ 
simplifies to 
\begin{equation} \label{densitybetahalf}
 f_{X_{\tp, \frac{1}{2}}}(x) = p f_1(x) + (1-p) f_2(x), 
\end{equation}
where 
\[
p= \frac{P_1 + P_2}{2} = \frac{1}{4} \left ( p_1+p_2+p_3+p_4\right ) = E[\tp].   
\] 
Therefore, the model written in \S \ref{1stmodelsubsec}, which was written 
intuitively, can be obtained as a special defuzzification of the fuzzy model treated 
in this subsection by setting $\beta = 1/2$. 
\end{remark}

Now we explain that the second model with a rectangular fuzzy number, 
written in \S \ref{2ndmodelsubsec}, can also be realized as a special case 
of the fuzzy model. 

\begin{remark} \label{model32relate}
If we assume that the fuzzy number 
$\tp$ in the fuzzy demand model $\tX$ given by 
\eqref{fuzzydemandintro} is rectangular as in the model 
in \S \ref{2ndmodelsubsec}, namely $\tp=$ $(p_m, p_m, p_M, p_M)$, 
then we have $P_1=2 p_m p_M$,  
$P_2$ $=$ $-2 p_m p_M+p_m+p_M$, and $P_3$ $=$ 
$2 \left(1-p_m\right) \left(1-p_M\right)$. Therefore, for the particular risk 
factor $\beta = 1/2$, the distribution of the defuzzification 
$X_{\tp, \frac{1}{2}}$ of $\tX$ is given by 
\[
 f_{X_{\tp, \frac{1}{2}}}(x) = p f_1(x) + (1-p) f_2(x), 
\]
where 
\[
p= \frac{P_1 + P_2}{2} = \frac{p_m + p_M}{2} = E[\tp]. 
\] 
Hence, the optimization in the second model written in 
\S \ref{2ndmodelsubsec} is also obtained as a special defuzzification 
of the fuzzy model by setting $\beta = 1/2$. 

\end{remark}

Using the above remarks, we can now explain how the three 
models in this article associated with the fuzzy weight $\tp$ 
reduce to the classical newsvendor model if there is no 
uncertainty in the choice of $\tp$. In fact, this happens 
when the risk factor $\beta = 1/2$ and the weight $\tp$ 
is taken to be an ordinary (or a {\it crisp}) number.   

\begin{remark}
If the fuzzy number $\tp$ in the fuzzy demand model $\tX$ given 
by  \eqref{fuzzydemandintro} is taken to be of the form 
$\tp = (p, p, p, p)$, where $p$ is in the 
interval $[0, 1]$, then $P_1 = 2 p^2$, 
$P_2 = -2 (p-1) p$, $P_3 = 2 (p-1)^2$, which yields $(P_1+P_2)/2=p$ 
and $(P_2+P_3)/2=1-p$. Then, if one takes 
the risk factor $\beta$ equal to $1/2$, considering Remarks 
\ref{model31relate} and \ref{model32relate}, the three fuzzy 
models considered in this article reduce to the classical 
newsvendor problem with the distribution function of the demand 
given by $f_X(x) = $ $p f_1(x) + (1-p) f_2(x)$.  
\end{remark}

Following the above remarks, only when we have $\beta = 1/2$ 
for the risk factor, the results that we have obtained by using the 
machinery of fuzzy random variables can be realized by the classical 
newsvendor analysis where the distribution function of the demand 
$X$ is of the form $p f_1(x) + (1-p) f_2(x)$.

Hence, a crucial result of our fuzzy analysis is the derivation of 
the distribution function $f_{X_{\tp, \beta}}$ given by 
\eqref{densitytpbeta} of the defuzzification $X_{\tp, \beta}$ 
of the fuzzy demand $\tX$. That is, in order to incorporate the 
uncertainties related to the weight of the online reviews in the 
formation of the demand, for a general risk factor $\beta$ in 
the interval $[0, 1]$, one can say that the optimal order quantity 
$Q^*_{\tp, \beta}$ is the solution of the newsvendor problem 
when the probability distribution function of the demand is taken to 
be the function $f_{X_{\tp, \beta}}$. Apparently, this function 
is by far more general than the distribution functions of the form 
$p f_1(x) + (1-p) f_2(x)$ that have been used classically and 
in the models written in \S \ref{1stmodelsubsec} and \S 
\ref{2ndmodelsubsec}. A great advantage of using the 
function $f_{X_{\tp, \beta}}$ is that it allows the newsvendor to 
use in the model the risk factor $\beta$ and the parameters of 
the fuzzy number $\tp$ that incorporates the range of online 
reviews and their weight in the formation of the demand. 
The point is that, as indicated by \eqref{densitybetahalf}, 
only when $\beta = 1/2$, the function $f_{X_{\tp, \beta}}$ is realized 
by the classical form $p f_1(x) + (1-p) f_2(x)$. Therefore, 
the nontrivial generalization $f_{X_{\tp, \beta}}$ that we have 
obtained by means of the fuzzy machinery is more powerful for 
the purpose of capturing the distribution of the demand.  
We shall provide in the sequel identification strategies for the 
risk factor $\beta$ and the parameters of $\tp$.   
 
\subsection{Risk-averse, risk-neutral, and risk-seeking decision making using the risk factor $\beta$}
Our objective in this subsection is to show that 
in our approach, in a natural way, we have 
introduced an objective function in which 
the parameter $\beta$ can encode the risk 
attitude of the decision maker. More importantly, for any chosen risk factor, we have an analytic 
formula for the optimizing order quantity. 
It is in fact crucial to incorporate 
risk attitude in studies related to the 
newsvendor problem, since it is often essential  for decision makers to have a control on 
their vulnerabilities to uncertainties. 
For studies in which mean-variance 
analysis is employed for the incorporation 
of risk attitude, we refer the reader to  
\cite{Choietal2008}, \cite{WuWangetal2009},  
and references therein.

In our approach, we have introduced a novel 
objective function introduced in \S \ref{betadefuzoptimizesubsec}, namely 
\begin{equation*} 
(1 - \beta)  \int_0^1 \bE[Y^1_\al] \, d\al+ \beta \int_0^1 \bE[Y^2_\al] \, d\al,   
\end{equation*}
which incorporates a risk factor $\beta$, and 
we have found analytically an order quantity 
$Q^*_{\tp, \beta}$ 
that optimizes our objective function for 
any given $\beta$. 
The following result shows that, using our model, the decision 
maker can increase its profit expectation 
by increasing the risk factor $\beta$. 
\begin{proposition}
\label{profitexpincreasinginbeta}
The profit expectation  
$\bE [ \pi (X_{\tp, \beta}) ]$ 
evaluated at the optimal value 
$Q^*_{\tp, \beta}$ introduces an increasing 
function of the risk factor $\beta$. 
\begin{proof}
Using \eqref{profitexpsimpleformula}, we can 
write: 
\[
\bE [ \pi (X_{\tp, \beta}) ] (Q^*_{\tp, \beta})
= 
(-M+V) \int_0^{Q^*_{\tp, \beta}} F_{X_{\tp, \beta}}(t) \, dt 
+ (M-C)Q^*_{\tp, \beta}.
\]
We find that 
\begin{eqnarray*}
&& \frac{d}{d\beta} \left (
\bE [ \pi (X_{\tp, \beta}) ] (Q^*_{\tp, \beta}) 
\right ) \\
&=& 
(V-M) \left(\int_0^{Q^*_{\tp, \beta}} 
\frac{\partial}{\partial \beta}F_{X_{\tp, \beta}}(t) \, dt+
\left ( \frac{d}{d\beta}Q^*_{\tp, \beta} \right )F_{X_{\tp, \beta}}(Q^*_{\tp, \beta})\right) \\
&+&(M-C) \frac{d}{d\beta}Q^*_{\tp, \beta} \\
&=& 
(V-M) \left(\int_0^{Q^*_{\tp, \beta}} \frac{\partial}{\partial \beta}F_{X_{\tp, \beta}}(t) \, dt \right )
\\
&=& (M-V) \left(\int_0^{Q^*_{\tp, \beta}} 
\left (J(t) - H(t) \right)\, dt \right ), 
\end{eqnarray*}
where we use the fact that 
$F_{X_{\tp, \beta}}(Q^*_{\tp, \beta}) = (M-C)/(M-V)$. 
Since $M-V >0$, and we know from Lemma 
\ref{GHJprop} that $H(t) \leq J(t)$ for all $t$, 
we can conclude that $\bE [ \pi (X_{\tp, \beta}) ] (Q^*_{\tp, \beta})$ is an increasing function 
of the risk factor $\beta$. 
\end{proof}
\end{proposition}

Based on Remark \ref{model31relate}, when 
$\beta = 1/2$, our fuzzy model is equivalent to a risk-neutral classical newsvendor problem. Hence, based on the results presented in Propositions \ref{Qincreasinginbeta} and \ref{profitexpincreasinginbeta}, a decision maker can use our fuzzy model 
as follows. In the risk-averse case, one can 
choose to have $0 \leq \beta < 1/2$, in the 
risk-neutral case $\beta = 1/2$, and for risk-seeking $ 1/2 < \beta \leq 1$.

\begin{example}
In the most risk-seeking case, the decision maker chooses $\beta = 1$, and models the
distribution of the demand with 
\[
f_{X_{\tp, 1}}(x)
=
P_1 F_1(x) f_1(x)+P_3 F_2(x) f_2(x)+P_2 \left(F_2(x) f_1(x)+F_1(x) f_2(x)\right).
\]
In this case, the order quantity is set by  
\begin{equation} \label{optimalfuzzybetaequalone}
Q^*_{\tp, 1} =  F^{-1}_{\tp, 1} \left ( \frac{M-C}{M-V} \right), 
\end{equation}
where
\[ 
F_{\tp, 1}(x) 
= 
\frac{1}{2} \left(P_1 F_1(x){}^2+2 P_2 F_1(x) F_2(x)+P_3 F_2(x){}^2\right). 
\]
\end{example}

\section{Example of application of the fuzzy GMM}
\label{numerical}
The aim of this section is twofold:
\begin{itemize}
\item To provide a concrete application scenario on how inventory managers could use both subjective and probabilistic measure of the demand uncertainty.
\item  To illustrate the implications of the mixture between the quantitative and subjective uncertainty of demand on inventory performance and to provide insights on the parameter inputs.

\end{itemize}
The source code of our simulations is available on  GitHub.\footnote{https://github.com/xxx/Fuzzy-Newsvendor-paper}

\subsection{Rationals for a fuzzy model}
We present a heuristics and an experimental study to legitimate the use of fuzzy numbers in our models, and give recommendations on how to derive them experimentally.

In Section \ref{fuzzydistsec}, we modelled the weight between online reviews and historical data in the formation of demand by a fuzzy number $\tp=(p_1, p_2, p_3, p_4)$.  We explain here how to derive $\tp$ from online data.

Consider a vendor's website, presenting its products for order and its customers online ratings for each of them, say on a 0 to 5 stars   scale. 
Consider a given product presented online, and let $\mr$ denote its average online ratings as displayed on its page.

We make the assumption that a proxy for the vendors' clients are website visitors. 
We split visitors into “customers” (on records for previous purchases), and “prospects” (with no prior purchases).

The heuristics behind the following recommendation is the following:
\begin{itemize}
    \item The left leg of $\tp$, ($p_1$, $p_2$), is related to the minimum weight for reviews {\it vs} historical distributions. This concerns {\em customers}, who are less sensitive to reviews and more driven by their experience with past products they bought and might need to order regularly.
    \item The right leg of $\tp$, ($p_3$,$p_4$), is related to the maximum weight for reviews {\it vs} historical distributions. This concerns {\em prospects}, who are sensitive to reviews since they have no experience with the company's products.
\end{itemize}

We further split customers into: {\em review insensitive} (ri) costumers who visit the website to order a product they know/need, and {\em review sensitive} (rs) costumers who might change their order based on reviews they see.
Let us denote by $n^v$ the total number of visitors on that page, $n^c$ the number of customers, and $n^p$ the number of prospects. We have $n^v=n^c+n^p=n^{ric} + n^{rsc} + n^p$, where $n^{ric}$ (resp. $n^{rsc}$) is the number of ri-customers (resp. rs-costumers).

The buying inclination of an rs-customer or a prospect checking up the product on that page is modelled with a trapezoidal fuzzy number $\tq =(q_1, q_2, q_3, q_4)$, with integers $0 \leq q_1 \leq q_2 \leq q_3 = q_4 =5$ as follows:
\begin{itemize}
    \item The visitor is certain {\it not} to place an order if  $\mr \le q_1$.
    \item The visitor is uncertain if $q_1 < \mr < q_2$.
    \item The visitor is certain to place an order if $\mr \ge q_2$.
\end{itemize}

For $rs$-customers the associated fuzzy numbers are denoted $\tq^{rsc}_i=(q_{i1}^{rsc},q_{i2}^{rsc},5,5)$, $i=1,\dots, n^{rsc}$, and for prospects $\tq^p_i=(q_{i1}^p,q_{i2}^p,5,5)$, $i=1,\dots, n^p$.
We assume that the expectations of prospects are higher than those of rsc-customers, that they need extra incentive to place an order. 
It follows that we model the $q_{ij}^{rsc}$ and $q_{ij}^p$ by Gaussian independent random variables with means satisfying
\[
 0\le  \mathbb{E}(q_{.1}^{rsc}) < \mathbb{E}(q_{.2}^{rsc}) < \mathbb{E}(q_{.1}^{p}) < \mathbb{E}(q_{.2}^{p}) \le 5 \,,
\]
and the same variances.

We further split the $n^{rsc}=n^{rsc}_0+n^{rsc}_1+n^{rsc}_2$ $rs$-customers into $n^{rsc}_0$ who did not place an order, $n^{rsc}_1$ who ordered without hesitation ({\it i.e.} for which $q_{i2}^{rsc}\le \bar{m}$), and $n^{rsc}_2$ who hesitated before ordering ({\it i.e.} for which $q_{i1}^{rsc}\le \bar{m} < q_{i2}^{rsc} $). 
Also, we similarly split the $n^{p}=n^{p}_0+n^{p}_1+n^{p}_2$ prospects into $n^{p}_0$ who did not order, $n^{p}_1$ who ordered without hesitation ({\it i.e.} for which $q_{i2}^{p}\le \bar{m}$), and $n^{p}_2$ who hesitated before ordering ({\it i.e.} for which $q_{i1}^{p}\le \bar{m} < q_{i2}^p $).

In order to derive a fuzzy $\tp=(p_1,p_2, p_3,p_4)$, the first step is to define a crip number $p_0$ simply weighting reviews vs historical demand as the ratio:
\begin{equation}
\label{p0}
p_0 = 
\frac{\# \{ \text{\rm customers who ordered b/c of reviews} \}}
{ \# \{ \text{\rm customers who ordered}\} }
= \frac{ n_1^{rsc}+n^{rsc}_2 }
{ n^{ris} + n_1^{rsc}+n_2^{rsc}} \,.
\end{equation}
Note that $p_0$ is related to $\tp$ via the defuzzification condition:
\begin{equation}
\label{p0-defuz}
    p_0 = \frac{p_1+p_2+p_3+p_4}{4}\,.
\end{equation}

Following the above heuristics, for the left leg of $\tp$ we set:
\[
p_1 = \alpha\frac{ n_1^{rsc} }{n}\,, \quad 
p_2 = \alpha \frac{n_1^{rsc}+n_2^{rsc}}{n}\,,
\]
and for the  right leg of $\tp$:
\[
p_3 = p_2 + \alpha\frac{ n_1^{p}}{n} \,, \quad
p_4 = p_2 + \alpha\frac{ n_1^{p}+n_2^p }{n}\,,
\]
where $n = n^{ric}+n_1^{rsc}+n_2^{rsc}+n_1^{p}+n_2^p$ is the total number of visitors who placed an order, and \(\alpha\) is a scaling factor which is computed from \eqref{p0} and \eqref{p0-defuz}: \(\alpha = 4np_0/(4n_1^{rsc}+3n_2^{rsc}+2n_1^p+n_2^p)\).
One easily checks that \(0\le p_1\le p_2\le p_3 \le p_4\le 1\).

\begin{figure}
    \centering
    \subfloat[\centering Fuzzy number $\tp=(0.57, 0.67, 0.75, 0.86)$ simulated from $N=10000$ online visitors, with mean ratings $\mr=3.5$. ]{{\includegraphics[width=5cm]{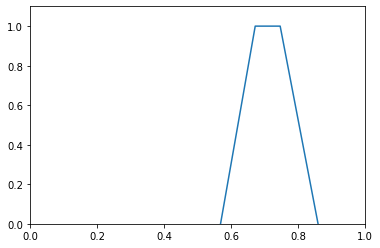} }}%
    \qquad
    \subfloat[\centering Simulated $\tp =(p_1, p_2, p_3, p_4)$ as a function of mean ratings $\mr$ (from 0.5 to 4.5 stars, with step 0.01) ]{{\includegraphics[width=5cm]{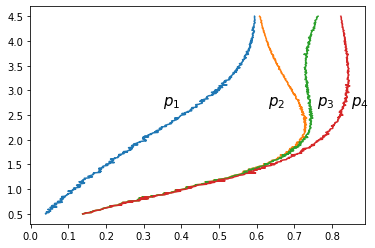} }}%
    \caption{Simulation to derive $\tp$}%
    \label{fig:simulation}
\end{figure}

We ran a numerical simulation, with 10,000 visitors, a prospect-customer ratio of 0.2, a ric-rsc ratio of 0.3, and Gaussian $q_{.i}$'s with means 1.5, 2.5 (rsc-customers), 3, 4 (prospects) and standard deviations 1, and mean ratings $\mr$ equal to 3.5 stars. See Figure~\ref{fig:simulation}(a).

We show in Figure \ref{fig:simulation}(b) how $\tp =(p_1, p_2, p_3, p_4)$ varies as a function of mean ratings $\mr$.

\subsection{Numerical analysis}
We consider now a numerical analysis to derive insights about the fuzzy GMM modelling of the newsvendor problem. For this purpose, we consider two main sources of demand information used by the decision maker when forecasting the demand:
\begin{itemize}
\item A normally distributed demand  resulting from a forecasting exercise based on past sales. Such a forecast exercise leads to a demand distribution with mean ${\mu }_2=100$ and standard deviation ${\sigma }_2=20$.  
\item A normally distributed demand adjusting the former by taking into account extra marketing actions and/or the impact of  past customer reviews or expert reviews on sales. This leads to a second distribution of demand with higher or lower moments compared to the first one. We assume that these extra actions potentially boost the demand and lead to a demand distribution with the mean ${\mu }_1=200$ and standard deviation ${\sigma }_1=30$.    
\end{itemize}

Mixing these two demand distributions is not a straightforward exercise in practise. Then, when a decision maker is facing these two alternatives for the demand distribution, s/he could choose only one by considering for instance the extreme cases by setting $p=0$ or $p=1$. Or s/he builds a mixture with a given weight for each distribution. The weight given to each alternative is by nature a subjective choice made by the decision maker since it depends on the behavior of customers. Indeed, estimating the impact of a marketing campaign on future demand or deriving the impact of past customer reviews on the behavior of future customers is not straightforward to model in a probabilistic manner. Following the customer behavior described in the previous section, we assume that the customers' sensitivity to review feedback (or marketing campaign) leads to a fuzzy distribution of  $\tp=(p_1, p_2, p_3, p_4)$. For the numerical illustration purpose let us consider two fuzzy values of $\tp$: $\tp_{case1}=(0.1, 0.2, 0.4, 0.4)$ and $\tp_{case2}=(0.6, 0.7, 0.9, 0.95)$.

Regarding the cost parameters, we consider two cost structures: $C=10$, $M=50$, $V=5$ ($C=10$, $M=12$, $V=5$) for high (low, respectively) margin product with a unit overage cost higher (lower, respectively) than the underage unit cost.

\begin{figure}
    \centering
    \subfloat[\centering for different values of $\tp$ ]{{\includegraphics[width=5cm]{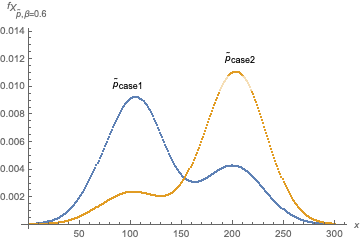} }}%
    \qquad
    \subfloat[\centering for different values of $\beta$ ]{{\includegraphics[width=5cm]{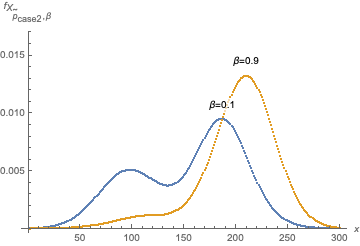} }}%
    \caption{Density function of the fuzzy GMM newsvendor}%
    \label{fig:density}%
\end{figure}

As proven in Section 4.4, we succeeded to show that the fuzzy GMM problem developed in this paper is equivalent to a classical newsvendor problem with a density function mixing the stochastic and fuzzy components of the demand as shown in 
 \eqref{densitytpbeta}. Given the mixture, the shape of this density function can be totally different from classical density functions usually used to solve the newsvendor problem. Indeed, the introduction of the GMM with fuzzy weights can lead to bi-modal demand distribution as illustrated in Figure \ref{fig:density}.

\begin{figure}
    \centering
    \subfloat[\centering Expected Demand]{{\includegraphics[width=5cm]{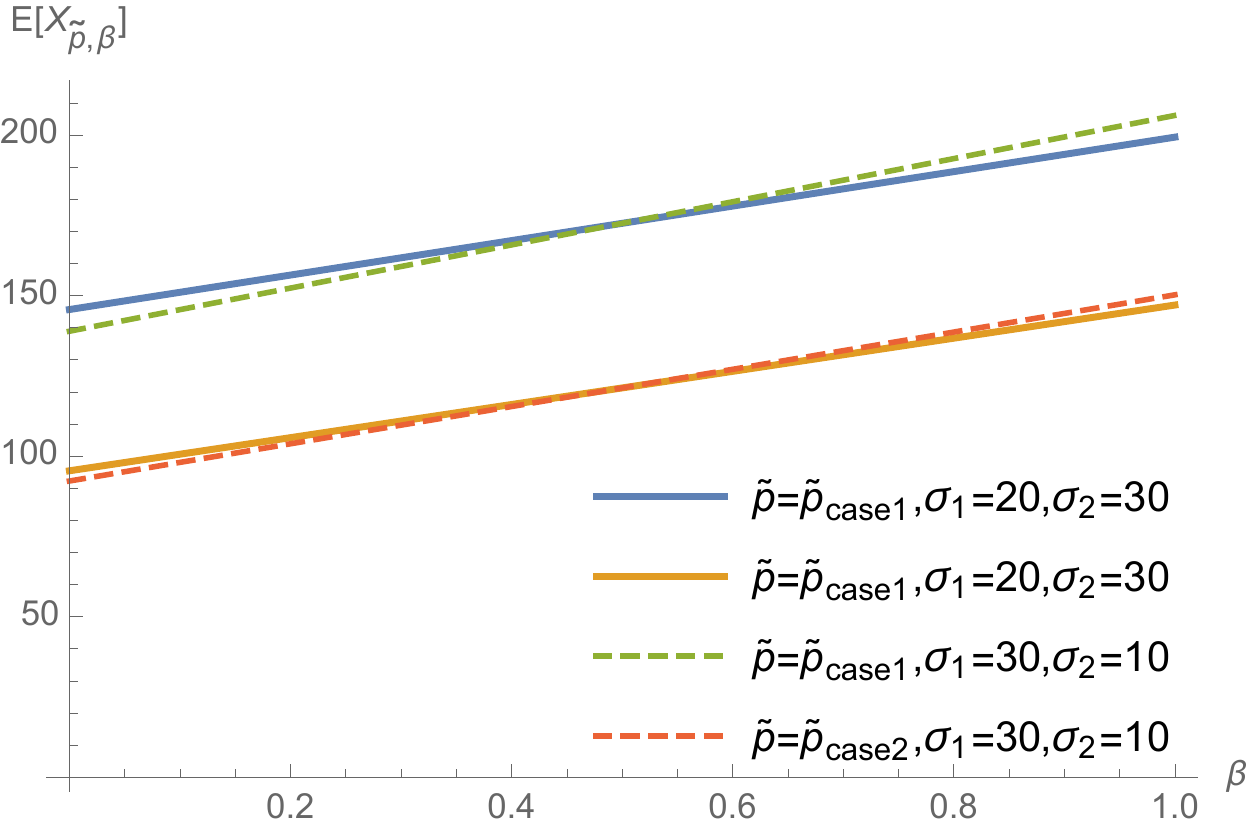} }}%
    \qquad
    \subfloat[\centering Demand Variance]{{\includegraphics[width=5cm]{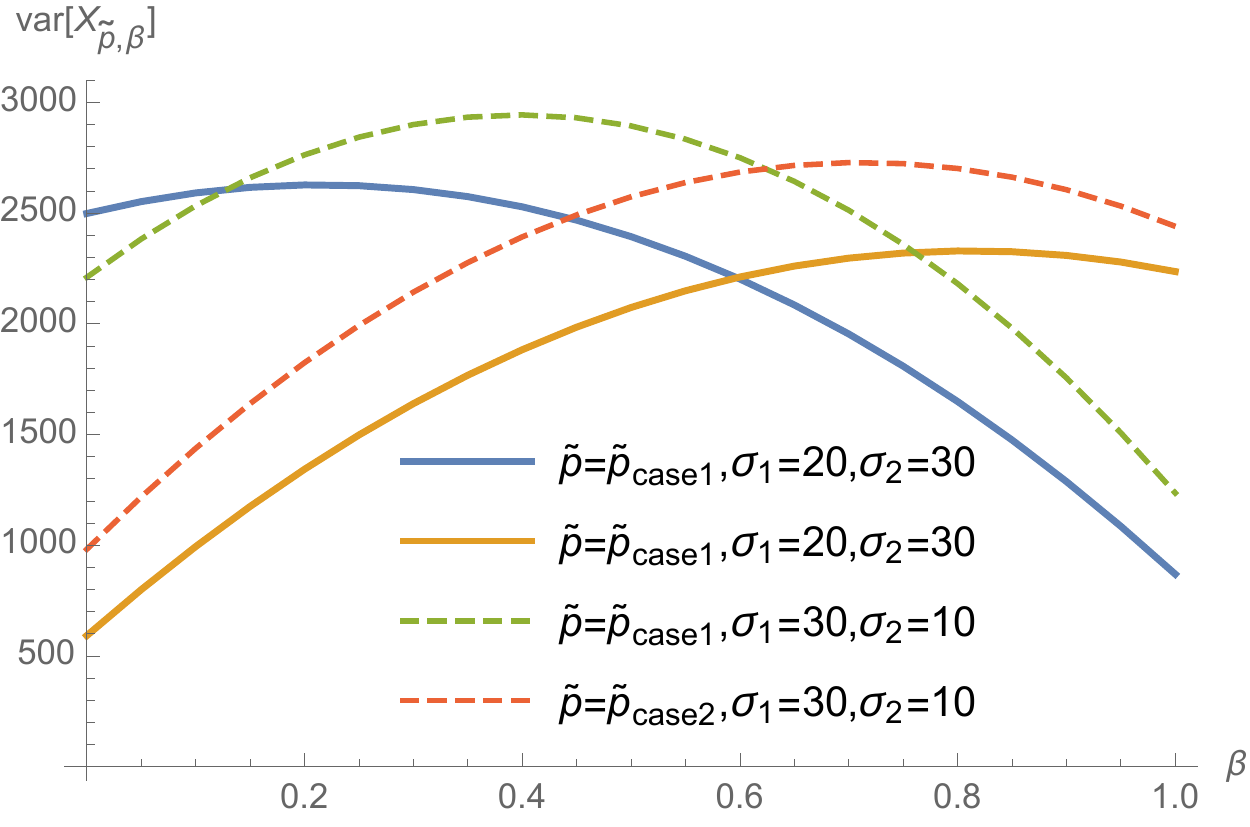} }}%
    \caption{Expectation and Variance of the fuzzy GMM demand }%
    \label{fig:expected}%
\end{figure}

It is worthwhile to notice that the higher the risk parameter $\beta$ is, the more important the second demand mode is when compared to the first mode. As a consequence, the value $\beta$ chosen by the decision maker impacts the order quantity and the associated profit expectation as well as the profit variance.
Indeed the choice made by the decision maker of the parameter $\beta$ drives the importance s/he gives to the right or the left side legs of the fuzzy number $\tp$. Indeed, the average $\bE[X_{\tp, \beta}]=\int_{-\infty}^{+\infty} xf_{X_{\tp, \beta}}(x) \, dx$  as well as the variance $\mathrm{Var}[X_{\tp, \beta}]=\int_{-\infty}^{+\infty} (x-\bE[X_{\tp, \beta}])^2 f_{X_{\tp, \beta}}(x) \, dx$  of the demand distribution $X_{\tp, \beta}$ are directly linked to $\beta$.
As illustrated in Figure \ref{fig:expected}, as intuitively expected, the expected demand resulting from the fuzzy GMM is increasing with $\beta$ since for higher $\beta$, the decision maker tends to give more importance to the right side leg of the distribution of $\tp$. Less intuitively, the variance of the fuzzy GMM demand is not necessarily in its lowest value for the risk neutral case ({\it i.e.} for  $\beta=0.5$) and might be lower in the extreme situations with risk-averse and risk-seeking choices of $\beta$. Obviously, the standard deviations of the probability distributions $\sigma_1$ and $\sigma_2$ impact directly the variance of the fuzzy GMM as illustrated in Figure \ref{fig:expected}(b), but they have a low impact on the expected values as depicted in Figure \ref{fig:expected}(a).

In order to highlight the added value of the fuzzy GMM modelling, one can assume that the decision maker might adopt three ways to handle the mixture of the two probabilistic demands weighted by the fuzzy number $\tp$: 
\begin{itemize}
    \item Ignoring one of the probability distributions by setting $p=0$ or $p=1$: these two cases model for instance a decision maker fully trusting or not the impact of past customer reviews on sales. Classical newsvendor solutions are then used here by considering as a demand density function either $f_1$ for $p=1$ or $f_2$ for $p=0$.
    \item Averaging the fuzzy number $\tp$ and then considering a GMM with a non fuzzy weight $p = E[\tp]$ (as detailed in Section \ref{1stmodelsubsec}). This case corresponds to the framework developed by \cite{rekik_enriching_2017} by setting the Bernoulli parameter of their model equal to  $E[\tp]$.
    \item Considering a fuzzy GMM with $\tp=(p_1, p_2, p_3, p_4)$, then choosing a risk factor $\beta$ and deriving the associated order quantity as developed in the present paper. 
\end{itemize}

For each case, the order quantity is calculated: $Q^*_{p=0}$ if $p=0$, $Q^*_{p=1}$ if $p=1$, $Q^*_{p = E[\tp]}$ if $p = E[\tp]$ and  $Q^*_{\tp, \beta}$ if $\tp=(p_1, p_2, p_3, p_4)$ for a given risk factor $\beta$. Then, the profit performance pertaining to each order quantity are derived by using the excepted profit and the profit variance of the fuzzy GMM $\bE [ \pi (X_{\tp, \beta}) ](\cdot)$,  $\mathrm{Var}[ \pi (X_{\tp, \beta}) ](\cdot)$, applied to each order quantity. 

Since $\bE [ \pi (X_{\tp, \beta}) ](\cdot)$ is optimized if the order quantity is set equal to $Q^*_{\tp, \beta}$, the three other cases, {\it i.e.} when the order quantities are equal to $Q^*_{p=j}$ $(j=0,1, E[\tp])$, which are then sub-optimal. We calculate their deviation from the optimal solution by using the following profit expectation ratios for $j=(0,1, E[\tp])$:

\begin{equation*} 
\textrm{Benefit of the fuzzy GMM} =  \frac{ \bE [ \pi (X_{\tp, \beta}) ](Q^*_{\tp, \beta})-\bE [ \pi (X_{\tp, \beta}) ](Q^*_{p=j}) }{\bE [ \pi (X_{\tp, \beta}) ](Q^*_{p=j}) }.
\end{equation*}

Similarly, for each $j=(0,1, E[\tp])$, we investigate the implication on the profit variance by using the following expression:
\begin{equation*} 
\textrm{Profit variance change} =  \frac{ \mathrm{Var}[ \pi (X_{\tp, \beta}) ](Q^*_{\tp, \beta})- \mathrm{Var}[ \pi (X_{\tp, \beta}) ](Q^*_{p=j}) }{\mathrm{Var}[ \pi (X_{\tp, \beta}) ](Q^*_{p=j}) }.
\end{equation*}

\begin{figure}
    \centering
    \subfloat[\centering $\tp=\tp_{case1}$]{{\includegraphics[width=5cm]{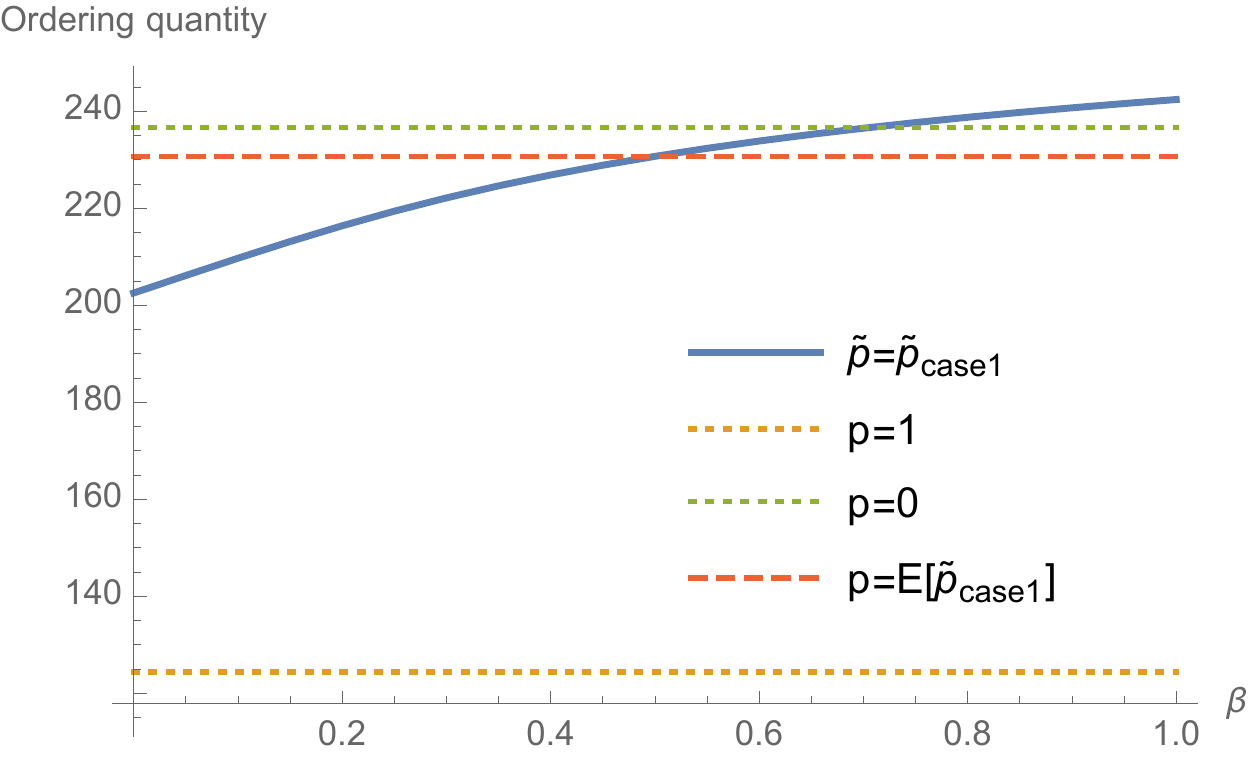} }}%
    \qquad
    \subfloat[\centering $\tp=\tp_{case2}$]{{\includegraphics[width=5cm]{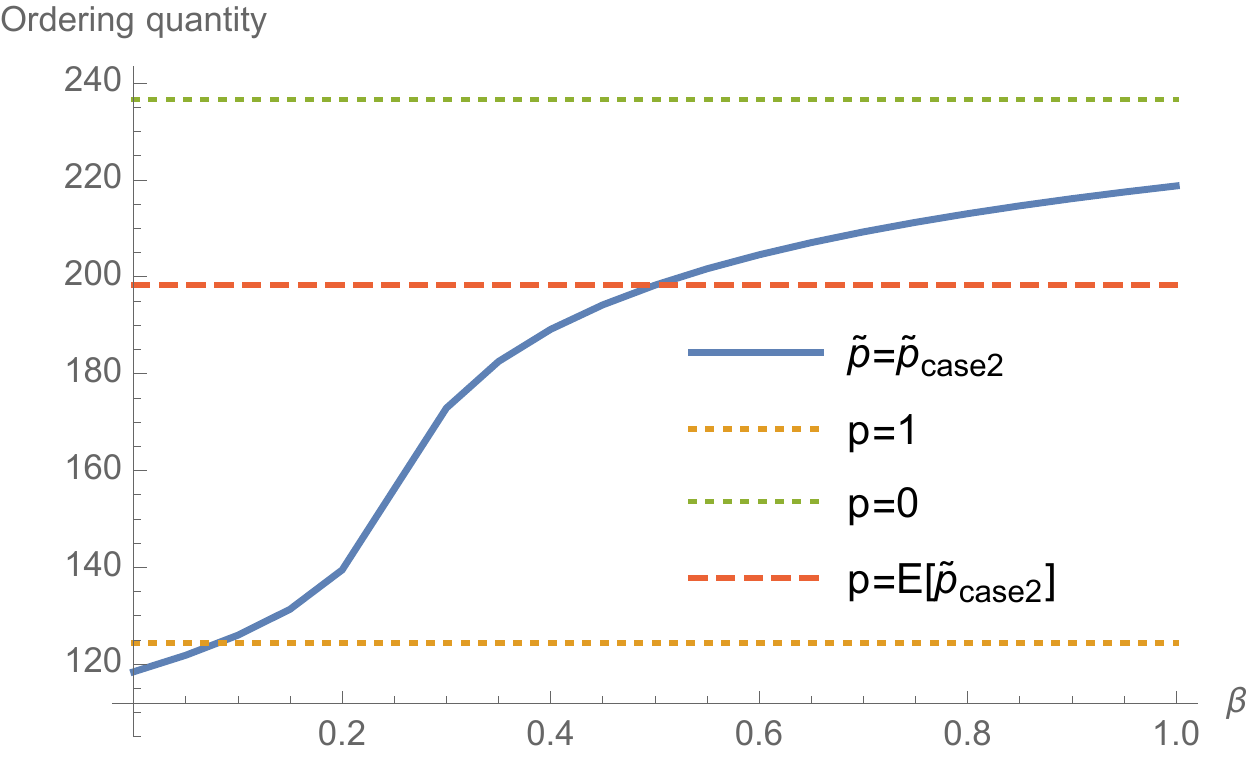} }}%
    \caption{Order quantity - high margin cost structure  }%
    \label{fig:8}%
\end{figure}

\begin{figure}
    \centering
    \subfloat[\centering $\tp=\tp_{case1}$]{{\includegraphics[width=5cm]{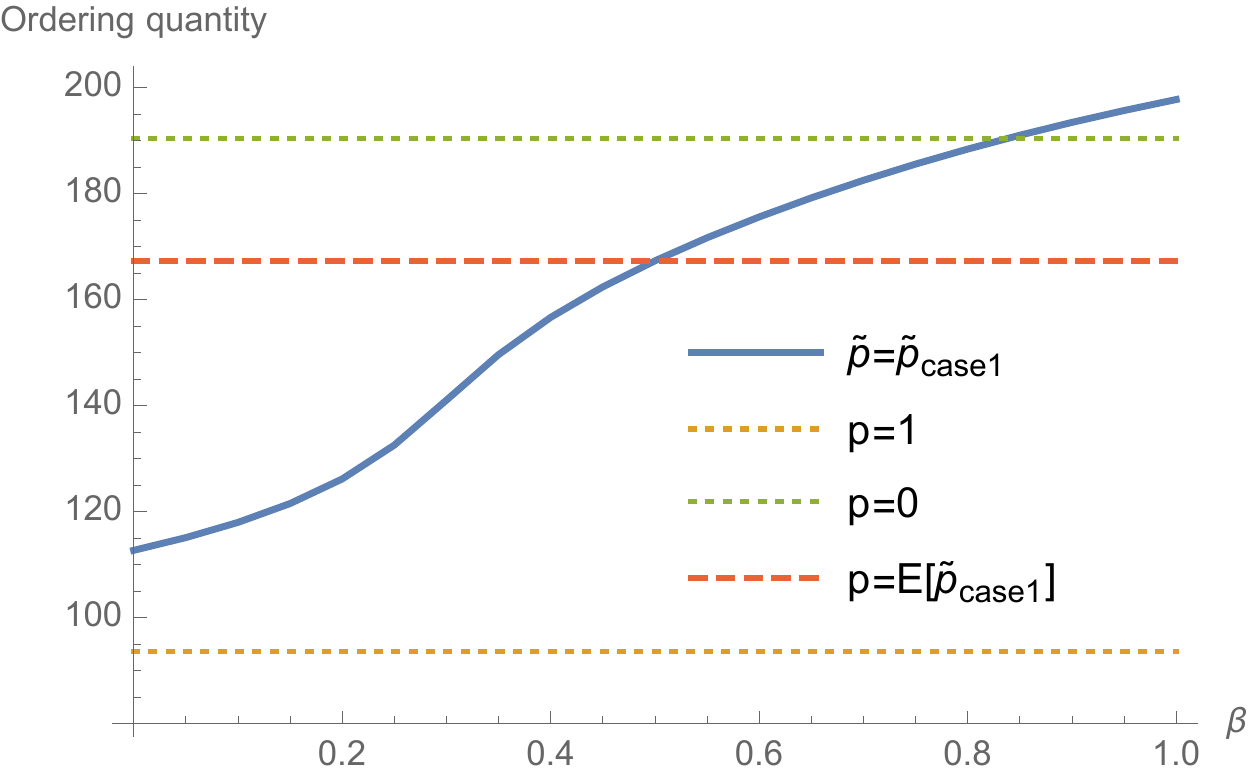} }}%
    \qquad
    \subfloat[\centering $\tp=\tp_{case2}$]{{\includegraphics[width=5cm]{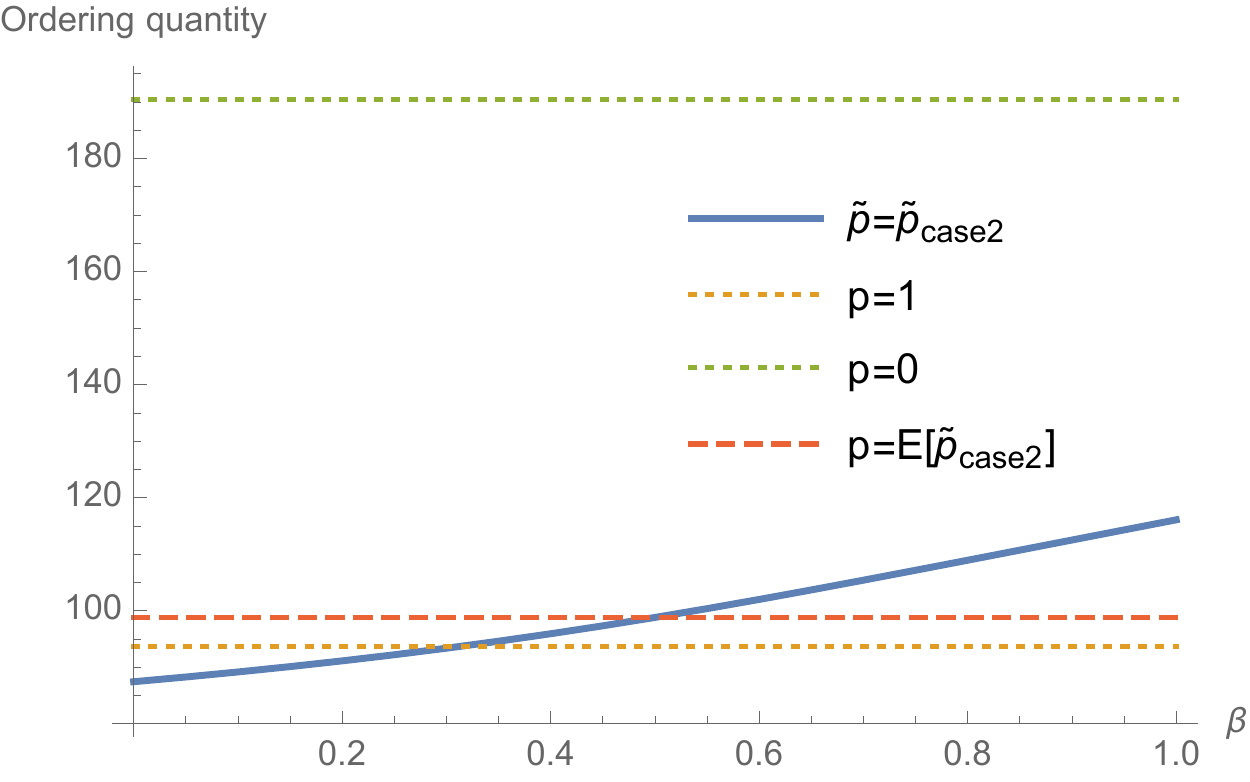} }}%
    \caption{Order quantity - low margin cost structure }%
    \label{fig:9}%
\end{figure}

Figures \ref{fig:8} and \ref{fig:9} compare the order quantities under the different values taken by $p$ for the high and low margin cost assumptions. Since $\bE[X_{\tp, \beta}]$ is increasing with $\beta$ as illustrated in Figure \ref{fig:expected}(a), the optimal order quantity of the fuzzy GMM, $Q^*_{\tp, \beta}$, increases with $\beta$ consequently. The slope of such increase is directly linked to the impact of $\beta$ on the demand variance $\mathrm{Var}[X_{\tp, \beta}]$ showed in Figure \ref{fig:expected}(b), but also to the cost structure of the product. While it is straightforward to analytically and numerically observe that the order quantity  $Q^*_{p = E[\tp]}$ is ranging between $Q^*_{p=0}$ and $Q^*_{p=1}$, this is not the case for $Q^*_{p = E[\tp]}$ which could fall outside this range for the extreme values of $\beta$. This is due to the fact that $Q^*_{p = \tp}$ involves not only the expectation of the fuzzy number $\tp$ but also its variability. It is known in the newsvendor literature that the order quantity increases (decreases) with the demand variance for high (low) margin cost structure. The graph pertaining to the case $p = E[\tp]$ crosses the fuzzy GMM one for  $\beta=0.5$: as discussed in Remark \ref{model32relate}, averaging the fuzzy number $\tp$ and using this average as a weight parameter is a particular case of our fuzzy GMM problem.

\begin{figure}
    \centering
    \subfloat[\centering Profit expectation]{{\includegraphics[width=5cm]{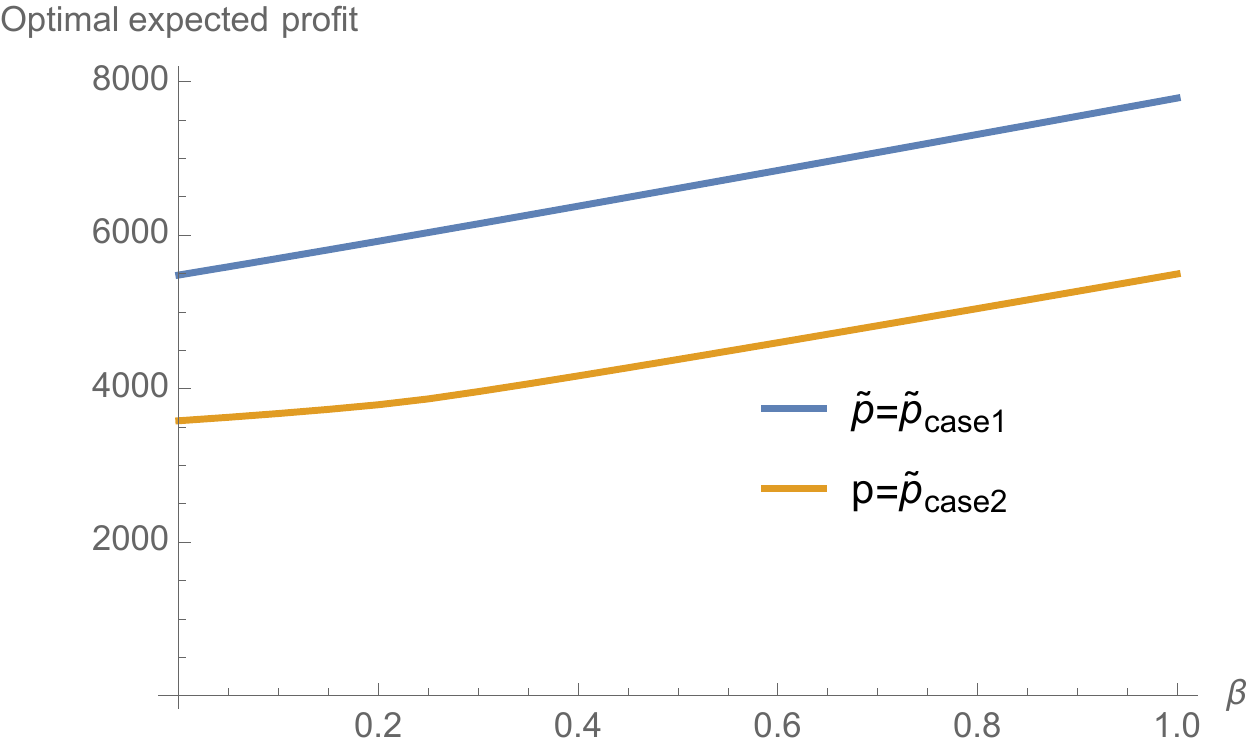} }}%
    \qquad
    \subfloat[\centering Profit variance]{{\includegraphics[width=5cm]{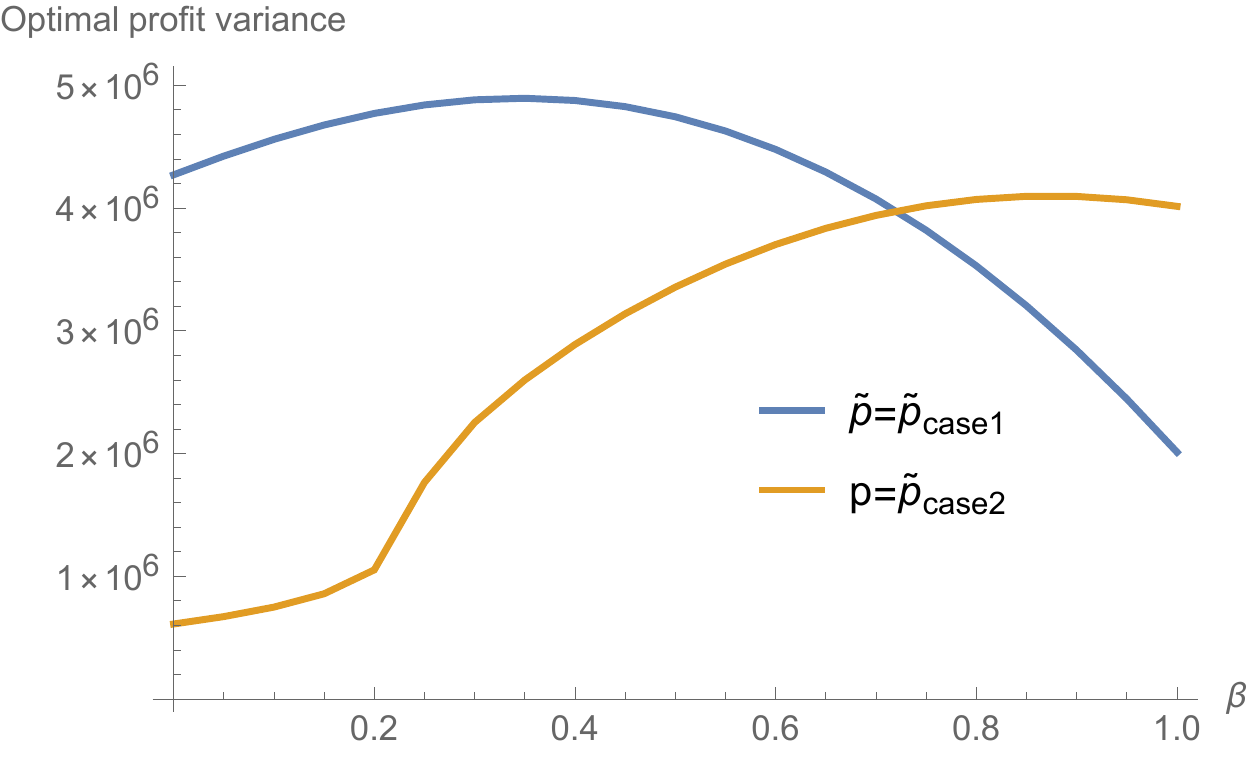} }}%
    \caption{Optimal Profit - high margin cost structure }%
    \label{fig:10}%
\end{figure}

\begin{figure}
    \centering
    \subfloat[\centering Profit expectation]{{\includegraphics[width=5cm]{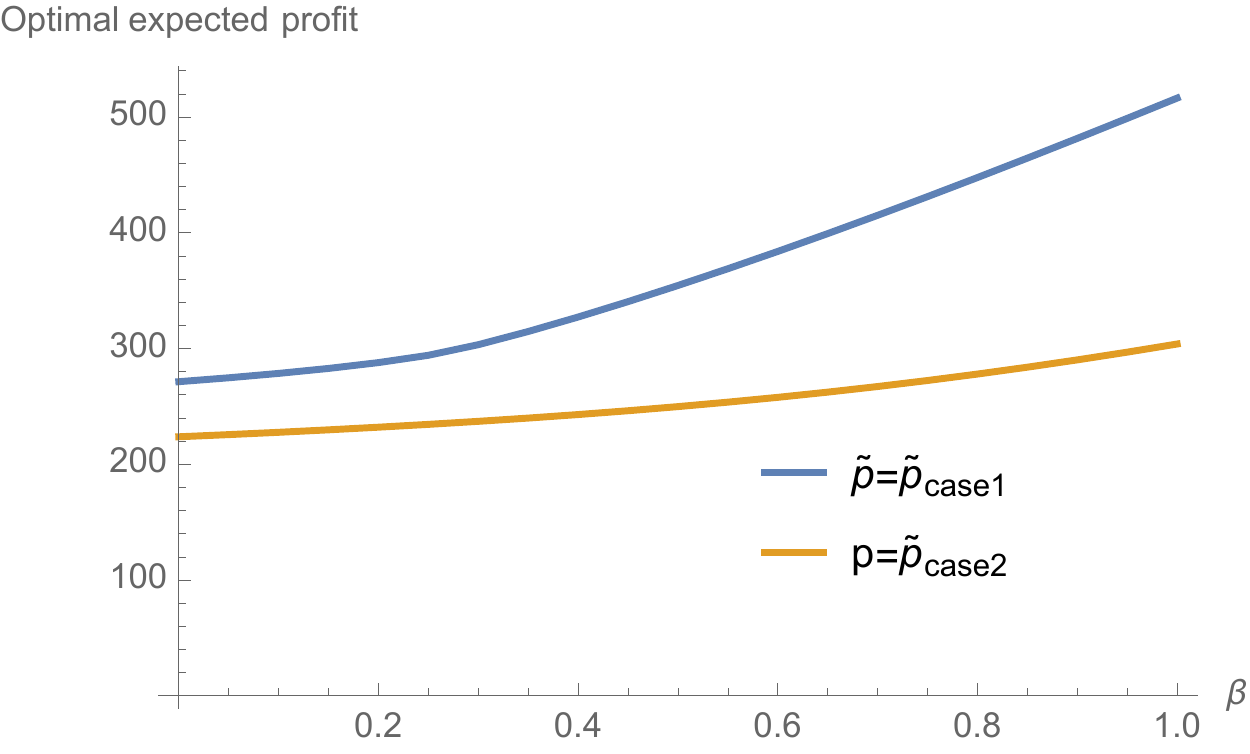} }}%
    \qquad
    \subfloat[\centering Profit variance]{{\includegraphics[width=5cm]{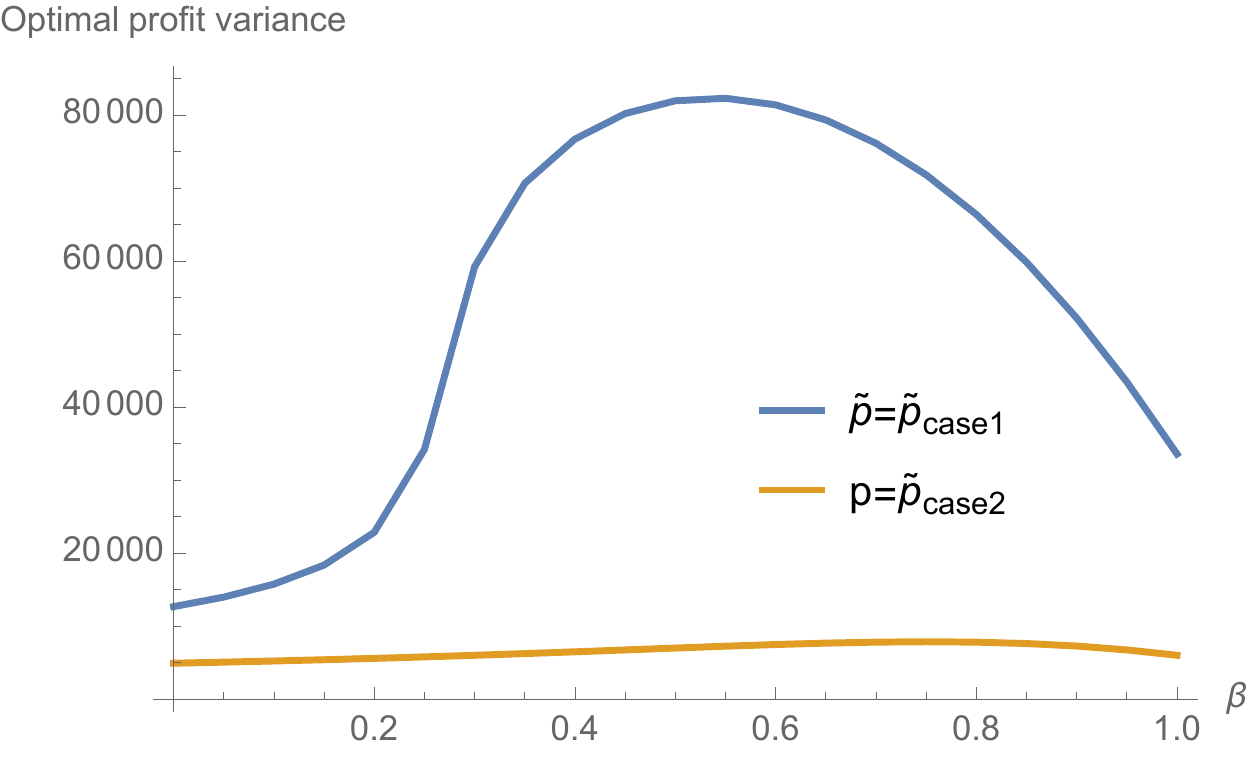} }}%
    \caption{Optimal Profit - low margin cost structure }%
    \label{fig:11}%
\end{figure}

By applying the fuzzy GMM for values of $\beta$ ranging from 0 to 1, the decision maker expects increasing profit for increasing $\beta$ as depicted in Figures \ref{fig:10}(a) and \ref{fig:11}(a): moving the cursor onto the right side of the fuzzy number increases the demand expectation and hence the profit expectation. Surprisingly, the profit variance (depicted in Figures \ref{fig:10}(b) and \ref{fig:11}(b)) is however not a monotonous function of  $\beta$ nor a convex function with a higher profit variance at the extreme values of $\beta$. In the contrast, it is worthwhile to notice that setting the $\beta$ around $0.5$, {\it i.e.} by giving similar importance to the left and right legs of the fuzzy number $\tp$, does not lead to a lower profit variance. For $\tp_{case1}=(0.1, 0.2, 0.4, 0.4)$, the profit variance is in its lowest value for $\beta=1$ under the high margin assumption and for $\beta=0$ under the low margin case and the result is reversed for $\tp_{case2}=(0.6, 0.7, 0.9, 0.95)$. By combining the outcome illustrated in the (a) and (b) sides of Figures \ref{fig:10} and \ref{fig:11}, the decision maker is able to choose the appropriate value of $\beta$ trading off its impact on the expectation and variance of the profit. Depending on the problem parameters, there are cases, such as $\tp_{case1}$ and high margin cost structure, where a choice of $\beta=1$ is beneficial for both the profit expectation and variance. There are also cases where a $\beta$ not in the extreme values is to be chosen depending on the risk sensitivity of the decision maker.

\begin{figure}
    \centering
    \subfloat[\centering Under case 2 $\tp_{case2}=(0.6, 0.7, 0.9, 0.9)$]{{\includegraphics[width=5cm]{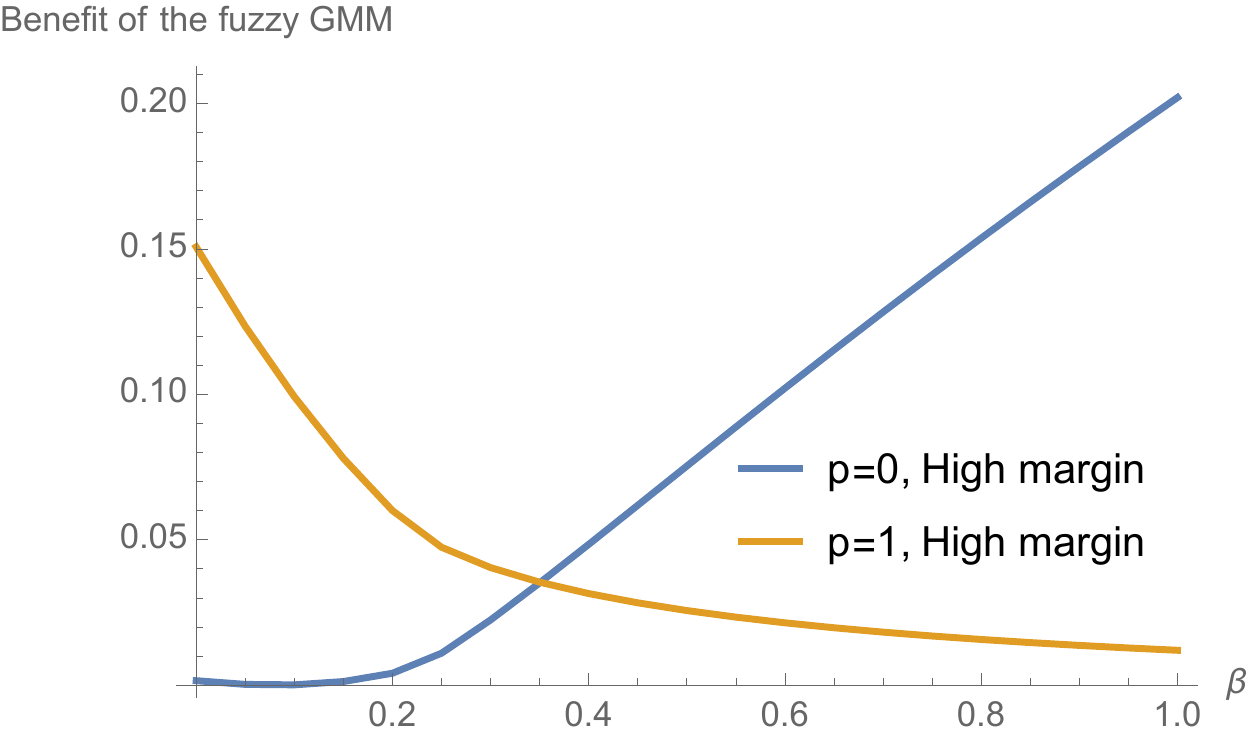} }}%
    \qquad
    \subfloat[\centering Under case 1 $\tp_{case2}=(0.1, 0.2, 0.4, 0.4)$]{{\includegraphics[width=5cm]{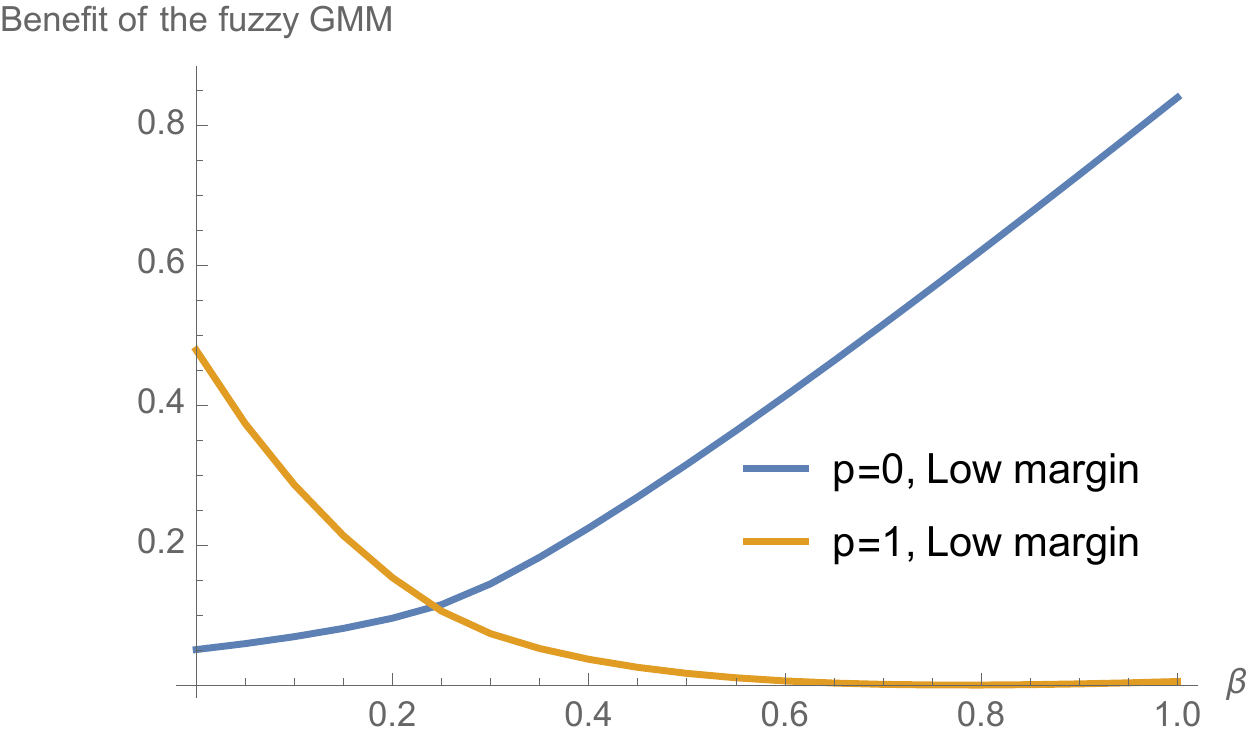} }}%
    \caption{Benefit of the fuzzy GMM if $p$ is set equal to 0 and 1  }%
    \label{fig:12}%
\end{figure}

As intuitively expected, Figure \ref{fig:12} shows that totally ignoring one of the probabilistic demand distributions leads to significant loss (up to $80\%$ for $\tp_{case1}$, low margin) compared to the case where the two distributions are mixed with a fuzzy number. This is a key message to practitioners who are facing increasing impact of social media, expert evaluation, and customer feedback on their statistical demand forecasts.

\begin{figure}
    \centering
    \subfloat[\centering Under the high margin cost]{{\includegraphics[width=5cm]{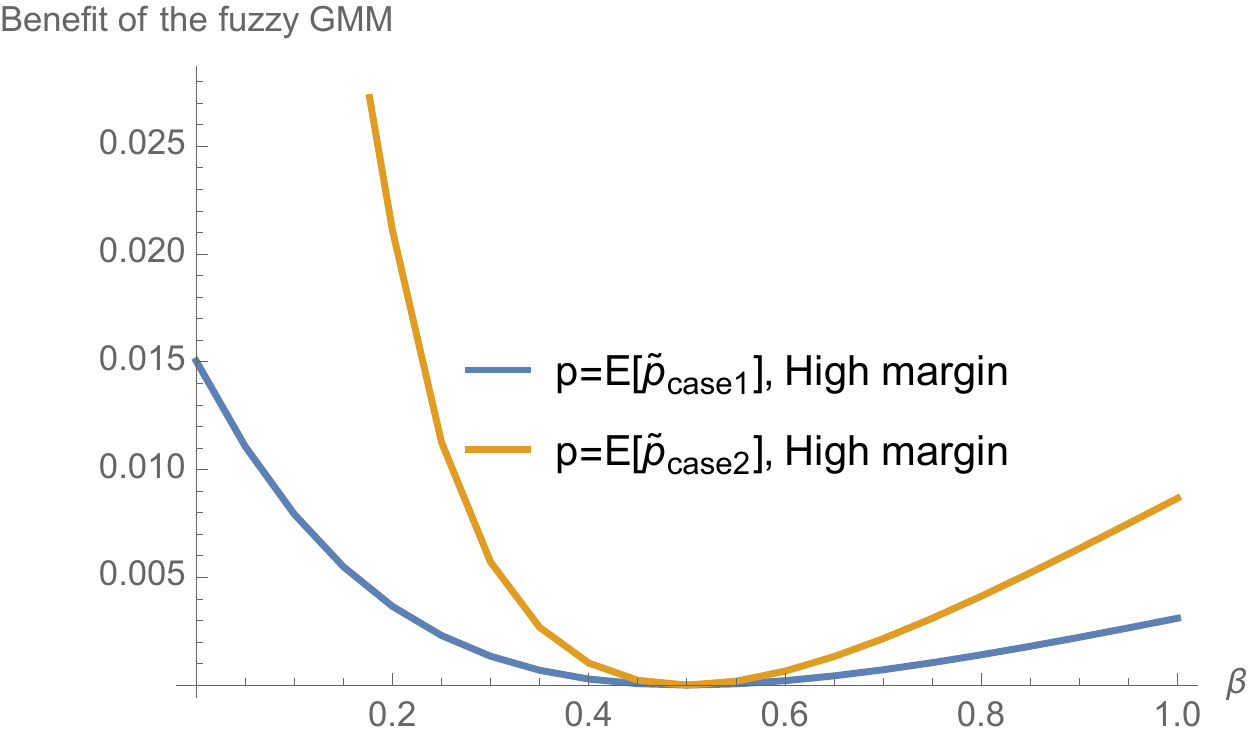} }}%
    \qquad
    \subfloat[\centering Under the low margin cost]{{\includegraphics[width=5cm]{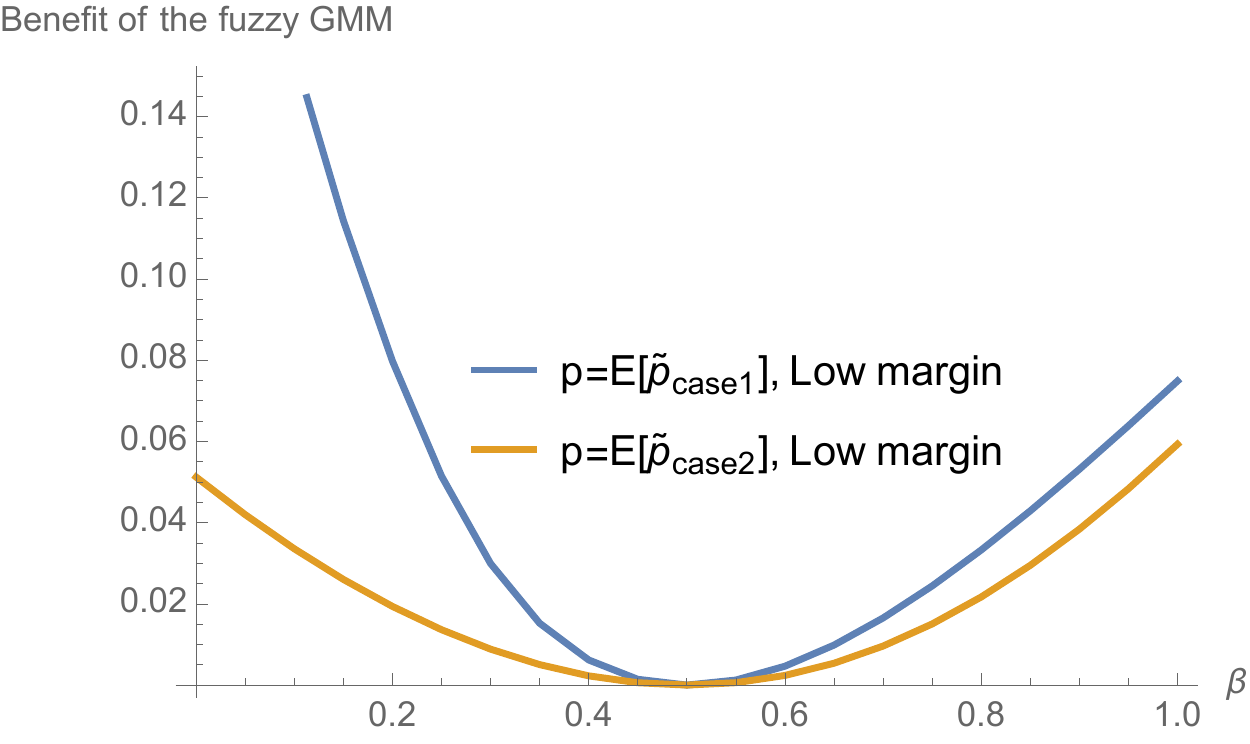} }}%
    \caption{Benefit of the fuzzy GMM if $p$ is set equal $p = E[\tp]$ }%
    \label{fig:13}%
\end{figure}

As mentioned above, the case where the fuzzy number $\tp$ is averaged and the average $E[\tp]$ is used as a weight of the GMM is a particular case of our fuzzy GMM problem. In order to measure the benefit related to the use of the fuzzy GMM instead of a GMM with a weight equal to the average of the fuzzy number, we compare the expected profit and profit variance in both cases. Since one is a particular case of the other, it is straightforward to observe a positive impact on the expected profit when using the fuzzy GMM (as illustrated in Figure \ref{fig:13}). More importantly, the benefit is quite meaningful particularly for the low margin cost assumption. Since, low margin products are sensitive to overstock penalties, the fuzzy GMM permits to fine tune the demand estimation with the factor $\beta$ and consequently to limit the overstock risk.

\begin{figure}
    \centering
    \subfloat[\centering Under case 2 $\tp_{case2}=(0.6, 0.7, 0.9, 0.9)$]{{\includegraphics[width=5cm]{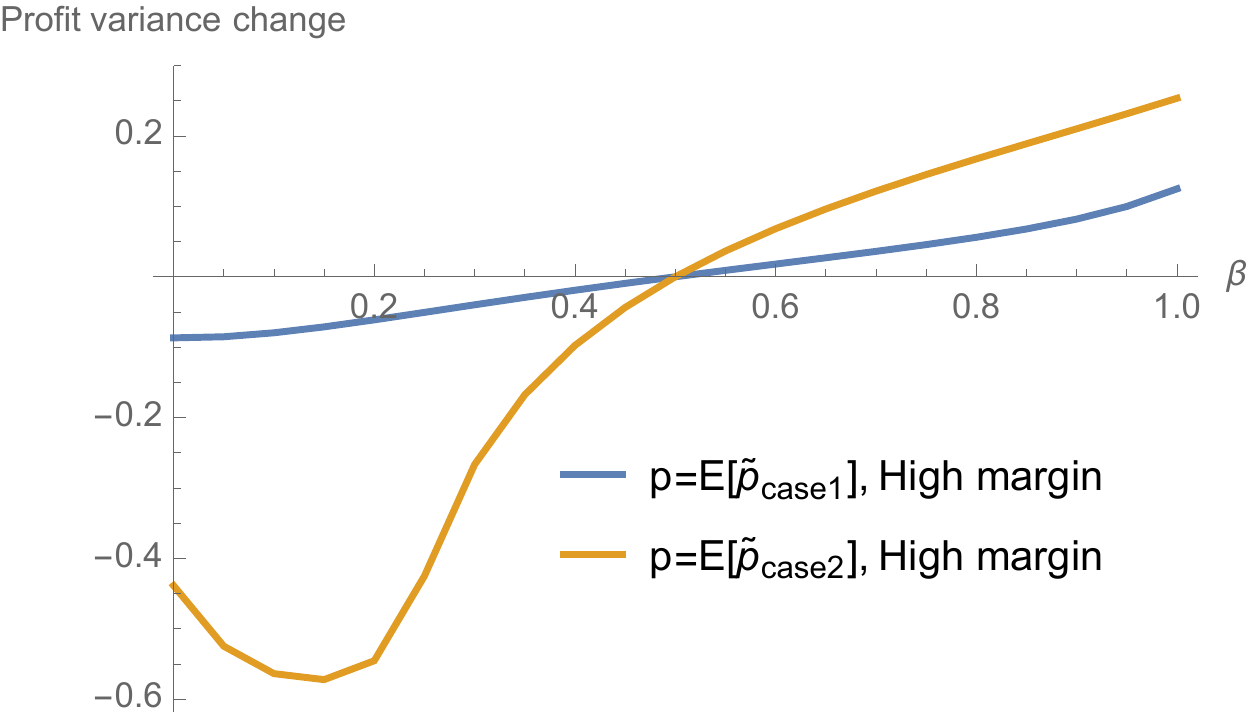} }}%
    \qquad
    \subfloat[\centering Under case 1 $\tp_{case1}=(0.1, 0.2, 0.4, 0.4)$]{{\includegraphics[width=5cm]{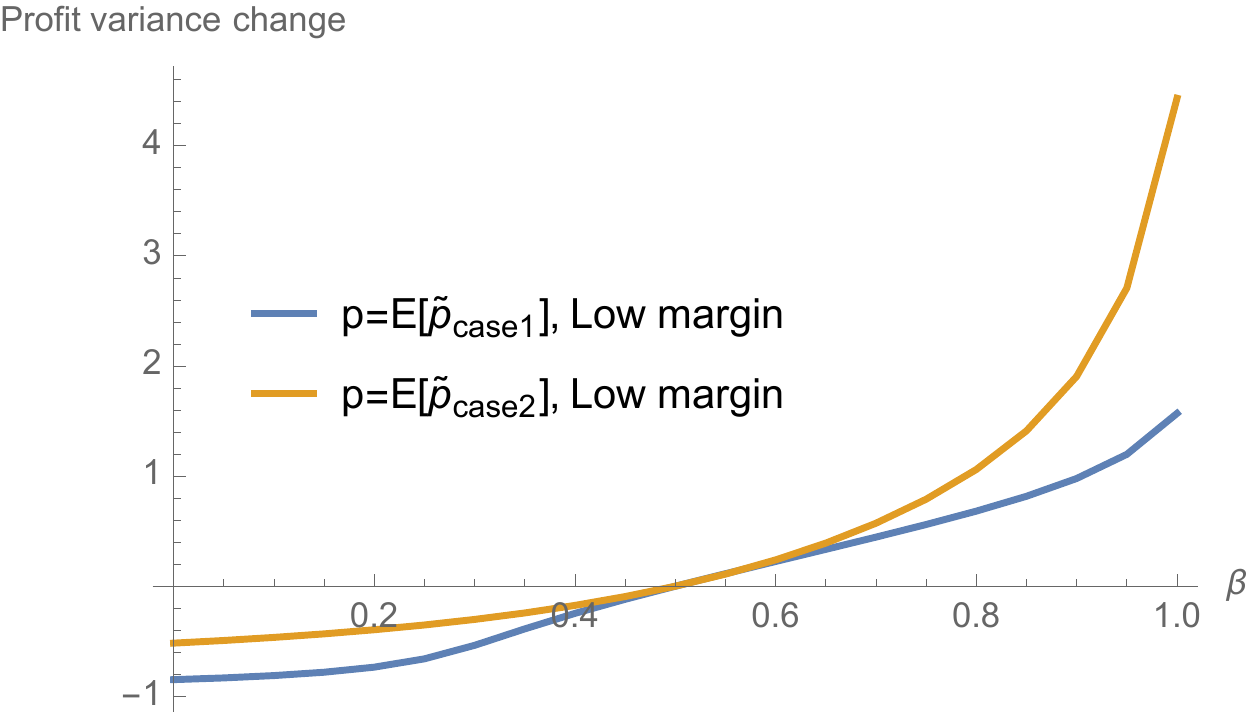} }}%
    \caption{Profit variance change of the fuzzy GMM if $p$ is set equal to $p = E[\tp]$ }%
    \label{fig:14}%
\end{figure}

Figure \ref{fig:14} shows that adopting the $\beta=0.5$ with an indifferent importance of the right and left legs of the fuzzy number $\tp$ can lead to a significant decrease of the profit variance under the risk-seeking profile of the decision maker. It is also worthwhile to notice a non-monotonous behavior of the profit variance impact. 

\section{Conclusion}
\label{conclusionssec}
Inspired by the increasing exposition of decision makers to both statistical and judgemental based sources of demand information, we developed in this paper a fuzzy GMM framework for the newsvendor permitting to mix probabilistic inputs with a subjective weight modelled as a fuzzy number. Thanks to a tractable mathematical application of the fuzzy machinery on the newsvendor problem, we derived the optimal ordering strategy taking into account both probabilistic and fuzzy components of the demand. We then reverse engineered the result to show that the fuzzy GMM problem could be interpreted as a classical newsvendor problem with a modified density function involving all the parameters of the probabilistic density functions and the legs of the fuzzy number. That is, mixing the different demand components could be interpreted as a requirement to revise the inventory policy by revisiting the ordering strategy of the newsvendor problem or a requirement to revise the forecasting exercise and then applying the classical newsvendor problem on the adjusted forecast outcome. The numerical study of the paper fitly provided the rationale of the use of fuzzy weight to the GMM. Customers sensitivity to review feedback or experts opinions is typically a perfect case study where customers demand is impacted by a fuzzy weighted GMM. The second part of the numerical study showed the applicability of the developed model and pointed out its added value compared to practices where the weights are either ignored or set equal to an average value avoiding then to integrate its variability in the inventory policy. The developed model permitted to show unusual multi-modal demand assumptions for the newsvendor problem and integrated the risk  measure on the profit.
As a future research agenda, one can target the extension of the model to a multi-period setting of the problem with a learning curve on the fuzzy number over periods. Indeed, the fuzzy weight was supposed to be known accurately in this paper which needs a learning exercise on the judgemental  based component of the demand.

\section*{Acknowledgements}
F.F. thanks the EMLYON Business School 
for their hospitality in November 2019 
and February 2020, and the great environment, 
where this work was partially carried out, 
and acknowledges the support from the Marie 
Curie/SER Cymru II Cofund
Research Fellowship 663830-SU-008. 



\end{document}